\documentclass[a4paper,11pt]{article}

\title{
\Large{Tessellation and Lyubich-Minsky laminations associated with quadratic maps, I: Pinching semiconjugacies
}
}
\author{
Tomoki Kawahira
\footnote{
Published in {\it Ergodic Theory \& Dynamical Systems}, 
{\bf 29} (2009), 579--612. This version contains a slight modification in Appendix A.2.
}
}

\date{\empty}
\usepackage{enumerate}
\usepackage{graphicx}
\usepackage{amsmath}
\usepackage{amscd}
\usepackage{amssymb}
\usepackage{amsfonts}
\usepackage{theorem}
\newtheorem{thm}{Theorem}[section]
\newtheorem{prop}[thm]{Proposition}
\newtheorem{lem}[thm]{Lemma}
\newtheorem{cor}[thm]{Corollary}
\theorembodyfont{\rmfamily}

\newtheorem{pf}{Proof.}

\newcommand{\thmref}[1]{Theorem \ref{#1}}
\newcommand{\propref}[1]{Proposition \ref{#1}}
\newcommand{\corref}[1]{Corollary \ref{#1}}
\newcommand{\lemref}[1]{Lemma \ref{#1}}

\newcommand{\figref}[1]{Figure \ref{#1}}
\newcommand{\secref}[1]{Section \ref{#1}}

\newcommand{\C}{\mathbb{C}}
\newcommand{\Cstar}{\mathbb{C}^\ast}

\newcommand{\Cbar}{\bar{\mathbb{C}}}
\newcommand{\Dbar}{\bar{\mathbb{D}}}
\newcommand{\R}{\mathbb{R}}
\newcommand{\D}{\mathbb{D}}
\newcommand{\Q}{\mathbb{Q}}
\newcommand{\Hyp}{\mathbb{H}}
\newcommand{\Z}{\mathbb{Z}}

\newcommand{\N}{\mathbb{N}}
\newcommand{\T}{\mathbb{T}}

\newcommand{\kakko}[1]{{\left( #1 \right)}}
\newcommand{\skakko}[1]{{\left\{ #1 \right\}}}

\newcommand{\llist}[1]{{{#1}_1, \ldots, {#1}_l}}

\newcommand{\QED}{\hfill $\blacksquare$}
\newcommand{\ee}{~=~}
\newcommand{\dee}{~:=~}

\newcommand{\lee}{~\le~}
\newcommand{\bs}[1]{\boldsymbol{#1}}
\newcommand{\parag}[1]{
\medskip
\noindent {\bfseries #1}
}
\newcommand{\card}{\mathrm{card}}

\newcommand{\rp}{\mathrm{Re}\,}
\newcommand{\ip}{\mathrm{Im}\,}
\newcommand{\diam}{\mathrm{diam}\,}
\newcommand{\depth}{\mathrm{depth}}
\newcommand{\val}{\mathrm{val}}

\newcommand{\Tess}{\mathrm{Tess}}

\newcommand{\id}{\mathrm{id}}

\newcommand{\al}{\alpha}
\newcommand{\gam}{\gamma}
\newcommand{\lam}{\lambda}
\newcommand{\Lam}{\varLambda}
\newcommand{\s}{\sigma}
\newcommand{\e}{\epsilon}
\newcommand{\cc}{\circ}
\newcommand{\fe}{f_\epsilon}

\newcommand{\fc}{f_c}

\newcommand{\phie}{\phi_\epsilon}
\newcommand{\He}{H_\epsilon}
\newcommand{\Kfc}{K_f^\circ}
\newcommand{\Kgc}{K_g^\circ}

\newcommand{\lbar}{{\bar{l}}}

\newcommand{\Tf}{T_f}

\newcommand{\Boe}{B}
\newcommand{\taue}{\tau_\epsilon}

\newcommand{\ue}{u_\epsilon}

\begin{document}

\maketitle

\begin{abstract}
We introduce \textit{tessellation} of the filled Julia sets for hyperbolic and parabolic quadratic maps. Then the dynamics inside their Julia sets are organized by tiles which work like external rays outside. We also construct continuous families of pinching semiconjugacies associated with hyperbolic-to-parabolic degenerations \textit{without} using quasiconformal deformation. Instead we use tessellation and investigation on the hyperbolic-to-parabolic degeneration of linearizing coordinates inside the Julia sets.
\end{abstract}

\section{Introduction}\label{sec_01}
After the works by Douady and Hubbard, dynamics of quadratic map $f=f_c:z \mapsto z^2+c$ with an attracting or parabolic cycle has been investigated in detail, because such parameters $c$ of $f_c$ are contained in the Mandelbrot set and they are very important elements that determine the topology of the Mandelbrot set. (See \cite{DH} or \cite{Mi}.)

The aim of this paper is to provide a new method to describe combinatorial changes of dynamics when the parameter $c$ moves from one hyperbolic component to another via a ``parabolic parameter" (i.e., $c$ of $f_c$ with a parabolic cycle). 
For example, the simplest case is the motion in the Mandelbrot set along a path joining $c=0$ and the center $c_{p/q}$ of the $p/q$-satellite component of the main cardioid via the root of the $p/q$-limb. 
In particular, we join them by the two segments characterized as follows:
\begin{itemize}
\item[(s1)] $c$ of $\fc$ which has a fixed point of multiplier $re^{2 \pi i p/q}$ with $0 < r \le 1$; and
\item[(s2)] $c$ of $\fc$ which has a $q$-periodic cycle of multiplier $1 \ge r > 0$. 
\end{itemize}
Note that we avoid the hyperbolic centers (i.e., $c$ of $f_c$ with superattracting cycle) because we regard them as non-generic special cases far away from parabolic bifurcations. 

\begin{figure}[htbp]
\centering{
\includegraphics[width=0.95\textwidth]{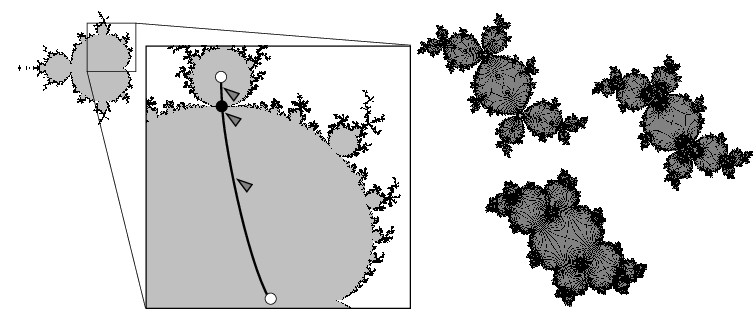}}
\caption{Chubby rabbits}\label{fig_rabbits}
\end{figure}

In the magnified box of \figref{fig_rabbits}, segments (s1) and (s2) for $p/q=1/3$ are drawn in the Mandelbrot set. By the Douady-Hubbard theory, the change of dynamics of $f=f_c$ on and outside the Julia set is described by external rays $R_f(\theta)$ with $\theta \in \T=\R/\Z$ and their landing points $\gam_f(\theta)$ satisfying $f(R_f(\theta))= R_f(2 \theta)$ and $f(\gam_f(\theta)) = \gam_f(2 \theta)$. For example, as $c$ moves from (s1) to (s2), the map $\gam_f:\T \to J_f$ loses the injectivity at a dense subset $\Theta_f$ of $\T$ consisting of the countably many angles that eventually land on $\skakko{1/7, 2/7, 4/7}$ by angle doubling $\delta: \theta \mapsto 2\theta$.

On the other hand, the dynamics inside the filled Julia set $K_f$ has no particular method to describe degeneration and bifurcation like external rays. However, as the pictures of filled Julia sets in \figref{fig_rabbits} (with equipotential curves drawn in) indicate, the interior of $K_f$ preserves a certain pattern along (s1) and (s2).

\parag{Degeneration pairs and tessellation.} 
In this paper, we introduce \textit{tessellation} of the interior $\Kfc$ of $K_f$ to detect hyperbolic-to-parabolic degeneration or parabolic-to-hyperbolic bifurcation of quadratic maps. 

Let $X$ be a hyperbolic component of the Mandelbrot set. By a theorem due to Douady and Hubbard \cite[Theorem 6.5]{Mi}, there exists the conformal map $\lam_X$ from $\D$ onto $X$ that parameterize the multiplier of the attracting cycle of $f=f_c$ for $c \in X$. Moreover, the map $\lam_X$ has the homeomorphic extension $\lam_X: \Dbar \to \bar{X}$ such that $\lam_X(e^{2 \pi i p/q})$ is a parabolic parameter for all $p, ~q \in \N$. A \textit{degeneration pair} $(f \to g)$ is a pair of hyperbolic $f=f_c$ and parabolic $g=f_\sigma$ where $(c, \sigma)=(\lam_X(re^{2 \pi i p/q}), \lam_X(e^{2 \pi i p/q}))$ for some $0<r<1$ and coprime $p, q \in \N$. By letting $r \to 1$, the map $f$ converges uniformly to $g$ on $\Cbar$ and we have a path which generalize segment (s1) or (s2). For a degeneration pair, we have the associated tessellations which have the same combinatorics:
\begin{thm}[Tessellation]\label{thm_1}
Let $(f \to g)$ be a degeneration pair. There exist families $\Tess(f)$ and $\Tess(g)$ of simply connected sets with the following properties:

\begin{itemize}
\item[\rm (1)] Each element of $\Tess(f)$ is called a {\rm tile} and identified by an {\rm angle} $\theta$ in $\Q/\Z$, a {\rm level} $m$ in $\Z$, and a {\rm signature} $\ast=+$ or $-$.
\item[\rm (2)] Let $T_f(\theta, m, \ast)$ be such a tile in $\Tess(f)$. Then $f(\Tf(\theta, m, \ast))=\Tf(2\theta, m+1, \ast)$.
\item[\rm (3)] The interiors of tiles in $\Tess(f)$ are disjoint topological disks. Tiles with the same signature are univalently mapped each other by a branch of $f^{-i} \cc f^j$ for some $i,~j > 0$;
\item[\rm (4)] Let $\Pi_f(\theta, \ast)$ denote the union of tiles with angle $\theta$ and signature $\ast$. Then its interior $\Pi_f(\theta, \ast)^\cc$ is also a topological disk and its boundary contains the landing point $\gam_f(\theta)$ of $R_f(\theta)$. In particular, $f(\Pi_f(\theta, \ast))=\Pi_f(2\theta, \ast)$.
\end{itemize}
The properties above holds if we replace $f$ by $g$. Moreover:
\begin{itemize}
\item[\rm (5)] There exists an $f$-invariant family $I_f$ of star-like graphs such that the union of tiles in $\Tess(f)$ is $\Kfc-I_f$. On the other hand, the union of tiles in $\Tess(g)$ is $\Kgc$.
\item[\rm (6)] The boundaries of $T_f(\theta, m, \ast)$ and $T_f(\theta', m', \ast')$ in $\Kfc -I_f$ intersect iff so do the boundaries of $T_g(\theta, m, \ast)$ and $T_g(\theta', m', \ast')$ in $\Kgc$.  
\end{itemize}
\end{thm}
\begin{figure}[htbp]
\centering{
\hspace*{0cm}\vspace*{0cm}
\includegraphics[width=\textwidth]{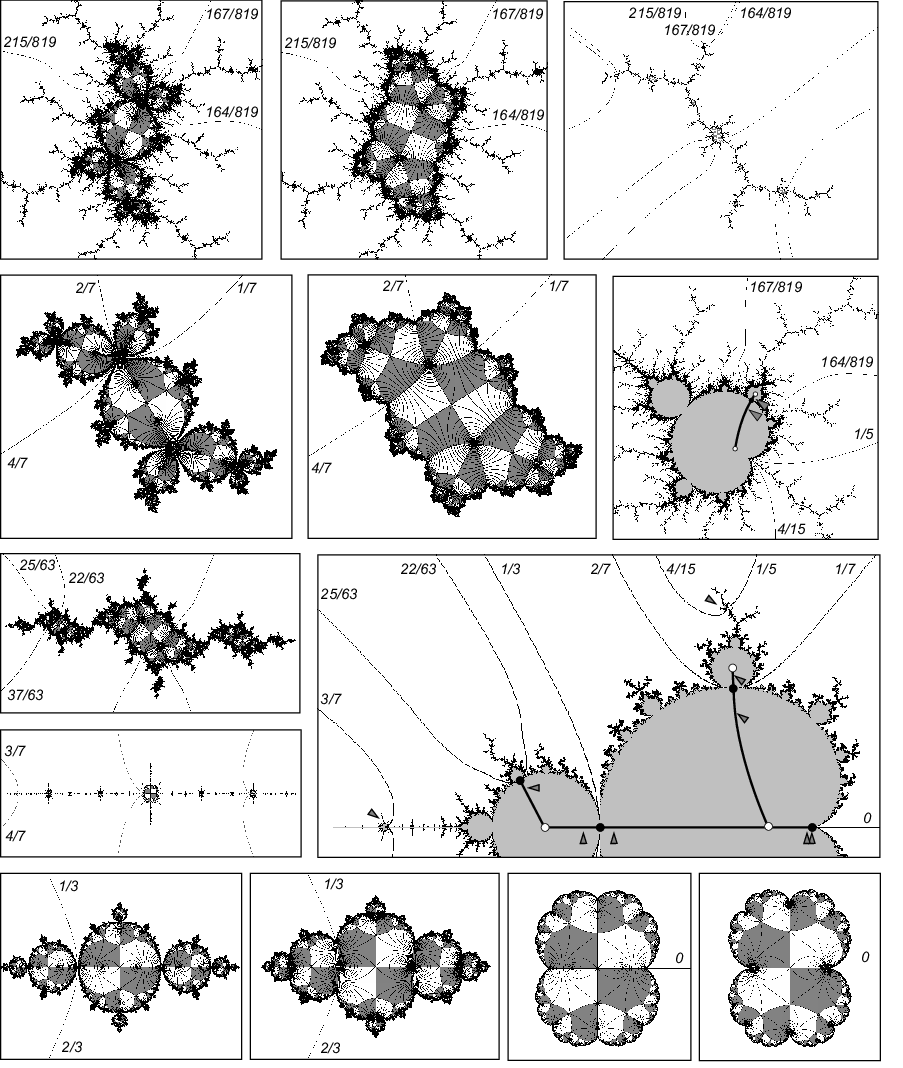}
}
\caption{Samples of tessellation. For the two figures at the upper left, parameters are taken from period 12 and 4 hyperbolic components of the Mandelbrot set as indicated in the figure of a small Mandelbrot set.}\label{fig_tessellation}
\end{figure}
Here angles of tiles must be the angles of external rays which eventually land on the parabolic cycle of $g$. (For example, if $(f \to g)$ are on (s1) or (s2) in \figref{fig_rabbits}, the set of angles of tiles coincides with $\Theta_f$.) See Sections \ref{sec_02} and \ref{sec_03} for construction of tessellation and \figref{fig_tessellation} for examples. One can find that the combinatorics of tessellations are preserved along (s1) and (s2). (This is justified in \secref{sec_04} more generally.) Since $f_c \in X-\skakko{\lam_X(0)}$ is structurally stable, we have the tessellation of $K_{f_c}$ with the same properties of $\Tess(f)$. 


\parag{Pinching semiconjugacy.}
As an application of tessellation, we show that there exists a pinching semiconjugacy from $f$ to $g$ for the degeneration pair $(f \to g)$. In Sections \ref{sec_04} and \ref{sec_05} we will establish:

\begin{thm}[Pinching semiconjugacy]\label{thm_2}
Let $(f \to g)$ be a degeneration pair. There exists a semiconjugacy $h: \Cbar \to \Cbar$ from $f$ to $g$ such that:
\begin{enumerate}[\rm (1)]
\item $h$ only pinches $I_f$ to the grand orbit of the parabolic cycle of $g$.
\item $h$ sends all possible $T_f(\theta, m, \ast)$ to $T_g(\theta, m, \ast)$, $R_f(\theta)$ to $R_g(\theta)$, and $\gam_f(\theta)$ to $\gam_g(\theta)$.
\item $h$ tends to the identity as $f$ tends to $g$. 
\end{enumerate}
\end{thm}
One may easily imagine the situation by seeing the figures of tessellation. As a corollary, we have convergence of the tiles when $f$ of $(f \to g)$ tends to $g$ (\corref{cor_conti_tiles}). 

We first prove the existence of $h$ with properties (1) and (2) in \secref{sec_04} (\thmref{thm_semiconj}) by using combinatorial properties of tessellation. Property (3) is proved in \secref{sec_05} (\thmref{thm_conti}) by means of the continuity results about the extended B\"ottcher coordinates (\thmref{thm_Boettcher_converges}) on and outside the Julia sets and the linearizing coordinates (i.e., the K\" onigs and Fatou coordinates) inside the Julia sets associated with $(f \to g)$ (\thmref{thm_linearization_converges}).

In Appendix, we will give some useful results on perturbation of parabolics used for the proofs. 

\parag{Notes.}
\begin{enumerate}
\item For any $f_c \in X-\skakko{\lam_X(0)}$, we have a semiconjugacy $h_c$ which has similar properties to (1) and (2) by structural stability. By results of Cui (\cite{Cui}), Ha\" {\i}ssinsky and Tan Lei(\cite{H3}, \cite{HT}), it is already known that such a semiconjugacy exists. Their method is based on the quasiconformal deformation theory and works even for some geometrically infinite rational maps. On the other hand, our method is faithful to the quadratic dynamics and the semiconjugacy is constructed in a more explicit way without using quasiconformal deformation. It is possible to extend our results to some class of higher degree polynomials or rational maps but it is out of our scope.

\item This paper is the first part of works on Lyubich-Minsky laminations. In \cite{LM}, they introduced the hyperbolic 3-laminations associated with rational maps as an analogue of the hyperbolic 3-manifolds associated with Kleinian groups. In the second part of this paper \cite{Ka4}, we will investigate combinatorial and topological change of 3-laminations associated with hyperbolic-to-parabolic degeneration of quadratic maps by means of tessellation and pinching semiconjugacies. 

\end{enumerate}

\parag{Acknowledgments.}
I would like to thank M. Lyubich for giving opportunities to visit SUNY at Stony Brook, University of Toronto, and the Fields Institute where a part of this work was being prepared. 
I am also grateful to the referee and Y.-C.~Chen for comments. 
This research is partially supported by JSPS Research Fellowships for Young Scientists, JSPS Grant-in-Aid for Young Scientists, Nagoya University, the Circle for the Promotion of Science and Engineering, and Inamori Foundation.

\section{Degeneration pair and degenerating arc system}\label{sec_02}
Segments (s1) and (s2) in the previous section are considered as hyperbolic-to-parabolic degeneration processes of two distinct directions. Degeneration pairs generalize all of such processes in the quadratic family. The aim of this section is to give a dichotomous classification of the degeneration pairs $\skakko{(f \to g)}$ and to define invariant families of star-like graphs (\textit{degenerating arc systems}) for each $f$ of $(f \to g)$.

\parag{Classification of degeneration pairs.}
We first fix some notation used throughout this paper. Let $p$ and $q$ be relatively prime positive integers, and set $\omega:=\exp(2 \pi i p/q)$. (We allow the case of $p=q=1$.) Take an $r$ from the interval $(0,1)$ and set $\lambda:=r\omega$. As in the previous section, we take a hyperbolic component $X$ of the Mandelbrot set. Then we have a \textit{degeneration pair} $(f \to g)$ that is a pair of hyperbolic $f=f_c$ and parabolic $g=f_\sigma$ where $(c, \sigma)=(\lam_X(re^{2 \pi i p/q}), \lam_X(e^{2 \pi i p/q}))$. 

For the degeneration pair $(f \to g)$, let $O_f:=\skakko{\llist{\al}}$ be the attracting cycle of $f$ with multiplier $\lam=r\omega$ and $f(\al_j)=\al_{j+1}$ (taking subscripts modulo $l$). Similarly, let $O_g:=\skakko{\beta_1, \ldots, \beta_{l'}}$ be the parabolic cycle of $g$ with $g(\beta_{j'})=\beta_{j'+1}$ (taking subscripts modulo $l'$). Let $\omega'= e^{2 \pi i p'/q'}$ denote the multiplier of $O_g$ with relatively prime positive integers $p'$ and $q'$. (Then $O_g$ is a parabolic cycle with $q'$ repelling petals.) 

Our fundamental classification is described by the following proposition:

\begin{prop}\label{prop_case_a_case_b}
Any degeneration pair $(f \to g)$ satisfies either
\begin{itemize}
\item[] {\rm {\bfseries Case (a)}:} ~~$q=q'$ and $l=l'$; or
\item[] {\rm {\bfseries Case (b)}:} ~~$q=1<q'$ and $l=l'q'$.  
\end{itemize}
For both cases, we have $lq=l'q'$.
\end{prop}
The proof is given by summing up results in sections 2, 4 and 6 of \cite{Mi}. For example, a degeneration pair $(f \to g)$ on segment (s1) (resp. (s2)) with $q>1$ is a Case (a) (resp. Case(b)) above. Degeneration pairs $(f_c \to f_\sigma)$ with $\sigma=1/4$ or $\sigma=-7/4$ satisfy $q=q'=1$ and thus Case (a).

\parag{Note on terminology.}
According to \cite{Mi}, a parabolic $g$ with $q'=1$ is called \textit{primitive}. The parabolic $g=f_\sigma$ with $\sigma=1/4$ is also called \textit{trivial}. For these $g$'s any degeneration pair $(f \to g)$ is automatically Case (a) by the proposition above. When we define tessellation for non-trivial primitive $(f \to g)$, we need an extra care.

\parag{Perturbation of $\bs{O_g}$ and degenerating arcs.}
For a degeneration pair $(f \to g)$ with $r \approx 1$, the parabolic cycle $O_g$ is approximated by an attracting or repelling cycle $O_f'$ with the same period $l'$ and multiplier $\lam' \approx e^{2 \pi i p'/q'}$. Let $\al_1'$ be the point in $O_f'$ with $\al_1' \to \beta_1$ as $r \to 1$. (cf. \cite[\S 4]{Mi})

In Case (a), the cycle $O_f'$ is attracting (thus $O_f'=O_f$) and there are $q'$ symmetrically arrayed repelling periodic points around $\al_1=\al_1'$. Then we will show that there exits an $f^{l'}$-invariant star-like graph $I(\al_1')$ that joins $\al_1'$ and the repelling periodic points by $q'$ arcs. In Case (b), the cycle $O_f'$ is repelling and there are $q'=l/l'$ symmetrically arrayed attracting periodic points around $\al_1'$. Then we will show that there exits an $f^{l'}$-invariant star-like graph $I(\al_1')$ that joins $\al_1'$ and the attracting periodic points by $q'$ arcs. 

In both cases, we define \textit{degenerating arc system} $I_f$ by 
$$
I_f \dee \bigcup_{n \ge 0} f^{-n}(I(\al_1')).
$$
The rest of this section is mainly devoted for the detailed construction of $I_f$, which has a role of parabolic cycle and its preimages. It would be helpful to see \figref{fig_deg_arc_sys} first, showing what we aim to have.


\subsection{External and internal landing}
First we consider the parameter $c$ on segment (s1) such that $f=f_c$ has an attracting fixed point $O_f=\skakko{\al_1}$ of multiplier $\lam=r \omega$, thus $c=\lam/2-\lam^2/4$. When $r$ tends to $1$, the hyperbolic map $f$ tends to a parabolic $g$ which has a parabolic fixed point $O_g=\skakko{\beta_1}$ with multiplier $\omega=e^{2 \pi i p/q}$. (Note that $q=q'$ and $l=l'(=1)$, thus Case (a) by \propref{prop_case_a_case_b}.) It is known that the Julia set $J_f$ of $f$ is a quasicircle, and the dynamics on $J_f$ is topologically the same as that of $f_0(z)=z^2$ on the unit circle. Since $J_f$ is locally connected, for any angle $\theta \in \R/\Z=\T$ its external ray $R_f(\theta)$ has a unique landing point $\gam_f(\theta)$. The same is true for $R_g(\theta)$, since $J_g$ is also locally connected. (See \cite[Expos\' e No.X]{DH}.)

\parag{External landing.}
By \cite[Theorem 18.11]{MiBook} due to Douady, there is at least one external ray with rational angle landing at $\beta_1$. Now \cite[Lemma 2.2]{GM} and the local dynamics of $\beta_1$ insure:

\begin{lem}\label{lem_external}
In the dynamics of $g$, there exist exactly $q$ external rays of angles $\theta_1, \ldots, \theta_q$ with $0 \le \theta_1< \cdots < \theta_q <1$ which land on $\beta_1$. Moreover, the map $g$ sends $R_g(\theta_j)$ onto $R_g(\theta_k)$ univalently (equivalently, $\theta_k = 2\theta_j$ modulo $1$) iff $k \equiv j+p \mod q$. 
\end{lem}
In particular, such angles are determined uniquely by the value $p/q \in \Q/\Z$. We take the subscripts of $\skakko{\theta_j}$ modulo $q$. For these angles, we call $p/q \in \Q/\Z$ the \textit{(combinatorial) rotation number}. Note that the external rays $\skakko{R_g(\theta_j)}$ divide the complex plane into $q$ open pieces, called \textit{sectors} based at $\beta_1$.

\parag{Internal landing lemma.} 
On the other hand, the set of landing points $\skakko{\gam_f(\theta_j)}$ of $\skakko{R_f(\theta_j)}$ is a repelling cycle of period $q$ and their corresponding rays do not divide the plane. However, they continuously extend and penetrate through the filled Julia set, and land at the attracting fixed point:

\begin{lem}[Internal landing] \label{lem_landing} 
For $\theta_1, \ldots, \theta_q$ as above, there exist $q$ open arcs \\
$I(\theta_1),\ldots ,I(\theta_q)$ such that:
\begin{itemize} 
\item For each $j$ modulo $q$, the arc $I(\theta_j)$ joins $\al_1$ and $\gam_f(\theta_j)$. 
\item $f$ maps $I(\theta_j)$ onto $I(\theta_k)$ univalently iff $k \equiv j+p \mod q$.
\end{itemize}
\end{lem}
There is no canonical choice for such arcs $\skakko{I(\theta_j)}$, but we will make a proper choice in the proof. 
An important fact is, the set $\skakko{I(\theta_j)\cup \gam_f(\theta_j) \cup R_f(\theta_j)}_{j=1}^q$ divide the plane into $q$ sectors based at $\al_1$. This is topologically the same situation as $g$. Indeed, as $r$ tends to $1$, we may consider that the arcs $\skakko{I(\theta_j)}$ constructed as above degenerate to the parabolic $\beta_1$.

\parag{Sketch of the proof.}
(See \cite[Lem. 2.3.3]{Ka3} for the detailed proof.) Let $w=\phi_f(z)$ be a linearizing coordinate near $\al_1$, where $f$ near $\al_1$ is viewed as $w \mapsto \lam w$. We can extend it to $\phi_f:\Kfc \to \C$ and normalize it so that $\phi_f(0)=1$ \cite[\S 8]{MiBook}. Now we pull-back the $q$th roots of the negative real axis in the $w$-plane, which are $q$ symmetrically arrayed invariant radial rays, to the original dynamics. Then we can show that the pulled-back arcs land at a unique repelling cycle with external angles determined by the rotation number $p/q$. In particular, they are disjoint from the critical orbit. \QED

\paragraph{Degenerating arcs.}
Note that in the construction of $\skakko{I(\theta_j)}$ above we make a particular choice of such arcs so that they are laid opposite to the critical orbit in the $w$-plane. We call these arcs \textit{degenerating arcs}.


\subsection{Degenerating arc system}
Let us now return to a general degeneration pair $(f \to g)$.

\parag{Renormalization.}
Let $B_1$ be the Fatou component containing the critical value $c$. We may assume that $B_1$ is the immediate basin of $\al_1$ for $f^l$. Then it is known that there exists a topological disk $U$ containing $B_1$ such that $f^l$ maps $U$ over itself properly by degree two. That is, the map $f^l: U \to f(U)$ is a quadratic-like map which is a renormalization of $f$. See \cite[\S 8]{Mi} or \cite[Partie 1]{H_Book}. In particular, the map $f^l: U \to f(U)$ is hybrid equivalent to $f_{1}(z)= z^2+c_1$ with $c_1=\lam/2-\lam^2/4$ in segment (s1), which we dealt with above. More precisely, the dynamics of $f^l$ near $B_1$ (resp. on $B_1$) is topologically (resp. conformally) identified as that of $f_{1}$ near $K_{f_1}$ (resp. on $K_{f_1}^\cc$).

\parag{Degenerating arc system.}
In $K_{f_1}$, we have $q$ degenerating arcs associated with the attracting fixed point of $f_1$. By pulling them back to the closure of $B_1$ with respect to the conformal identification above, we have $q$ open arcs $\skakko{I_j}_{j=1}^q$ which are cyclic under $f^l$. 

When $q=q'$ and $l=l'$, thus in Case (a), the arcs $\skakko{I_j}_{j=1}^q$ join $q'$ repelling points (cyclic under $f^l=f^{l'}$) and $\al_1=\al_1'$. In this case we define $I(\al_1')$ by the closure of the union of $\skakko{I_j}_{j=1}^q$. When $1=q<q'=l/l'$, thus in Case (b), we only have $I_1$ that joins the repelling point $\al_1'$ (fixed under $f^{l'}$) and $\al_1$. In this case we define $I(\al_1')$ by the closure of the union of $\skakko{f^{kl'}(I_1)}_{k=0}^{q'-1}$. In both cases, we have $I(\al_1')$ as desired. Now we define the \textit{degenerating arc system} of $f$ by
$$
I_f ~:=~ \bigcup_{n>0} f^{-n}(I(\al_1')).
$$
For $z \in I_f$, it is useful to denote the connected component of $I_f$ containing $z$ by $I(z)$. 

\begin{figure}[htbp]
\centering
\includegraphics[width=0.9\textwidth]{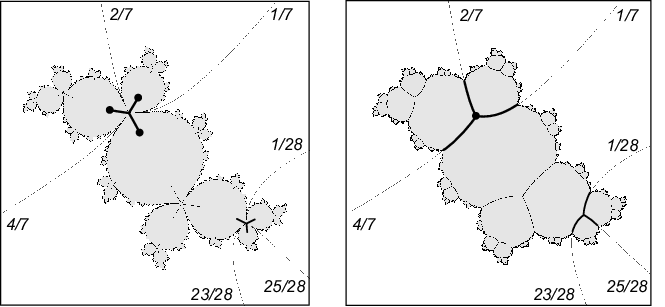}
\caption{Left, the Julia set of an $f$ in segment (s2) for $p/q=1/3$, and right, one in segment (s1), with their degenerating arc system roughly drawn in. Attracting cycles are shown in heavy dots. Degenerating arcs with types $\skakko{1/7, 2/7/ 4/7}$ and $\skakko{1/28, 23/28, 25/28}$ are emphasized.}\label{fig_deg_arc_sys}
\end{figure}

For later use, we define the set of all points that eventually land on the attracting cycle $O_f$, by $\al_f:= \bigcup_{n>0} f^{-n}(\al_1)$. Note that $I_f$ and $\al_f$ are forward and backward invariant, and disjoint from the critical orbit. In particular, for any $z \in I_f$, the components $I(z)$ and $I(\al_1')$ are homeomorphic. 
In Case (a), the points in $\al_f$ and the connected components of $I_f$ have one-to-one correspondence. In Case (b), however, they are $q'$-to-one correspondence. See \figref{fig_deg_arc_sys} and \propref{prop_valence}.

Correspondingly, for $g$ of degeneration pair $(f \to g)$ and one of its parabolic point $\beta_1 \in O_g$, we define 
$$
I_g ~:=~ \bigcup_{n>0} g^{-n}(\beta_1).
$$
We will see that this naturally corresponds to $I_f$ rather than $\al_f$.

\parag{Types.}
After \cite{GM}, we define the \textit{type} $\Theta(z)$ of $z$ in $J_f$ (or $J_g$) by the set of all angles of external rays which land on $z$. Let $\delta: \T \to \T$ be the angle doubling map. Since $J_f$ has no critical points, one can easily see that $\delta(\Theta(z))$ coincides with $\Theta(f(z))$. The same holds for $g$. 
By an unpublished work by Thurston, if $z$ in the quadratic Julia set does not have finite orbit, then $\card(\Theta(z)) \le 4$. See \cite{Ki} for a generalized statement and the proof. In our case, (pre)periodic points in $J_f$ and $J_g$ have uniformly bounded numbers of landing rays. Hence we have:
\begin{lem}\label{lem_type}
For any point $z$ in $J_f$ or $J_g$, the set $\Theta(z)$ consists of finitely many angles.
\end{lem}

We abuse the notation $\Theta(\cdot)$ like this: For any subset $E$ of the filled Julia set, its \textit{type} $\Theta(E)$ is the set of angles of the external rays that land on points in $E$. (So $\Theta(E)$ is empty when $E$ does not touch the Julia set.) For each $\zeta$ in $\al_f$, we formally define the type of $\zeta$ by $\Theta(\zeta):= \Theta(I(\zeta))$. Then one can easily see that $\delta^n(\Theta(\zeta))=\Theta(\al_1)$ for some $n>0$. We also set $\Theta_f:= \Theta(I_f)$ and $\Theta_g:= \Theta(I_g)$. We will show that $\Theta_f$ equals to $\Theta_g$ in the next proposition.

\parag{Valence.}
For any $\zeta \in \al_f$, the component $I(\zeta)$ of $I_f$ is univalently mapped onto $I(\al_1)$ by iteration of $f$. Thus the value $\val(f):=\card(\Theta(\zeta))$ is a constant for $f$. Similarly, since a small neighborhood of $\xi$ in $I_g$ is sent univalently over $\beta_1$ by iteration of $g$, the value $\val(g):=\card(\Theta(\xi))$ is constant for $g$. Now we claim:

\begin{prop}\label{prop_valence}
For any degeneration pair $(f \to g)$, we have $\Theta_f=\Theta_g$ and $\val(f)=\val(g)$. Moreover, 
\begin{itemize}
\item if $q=q'=1$ and $l=l'>1$ (thus Case (a) and non-trivial primitive), then $\val(g)=2$.
\item Otherwise $\val(g)=q'$.
\end{itemize}
\end{prop} 
We call $\val(f)=\val(g)$ the \textit{valence} of $(f \to g)$. Note that the valence depends only on $g$. 
The proof of this proposition uses some facts from \cite{Mi}.

\parag{Proof.}
The two possibilities of $\val(g)$ above are shown in \cite[Lemma 2.7, \S 6]{Mi}. If we show that $\Theta(\al_1')=\Theta(\beta_1)$, then $\Theta_f=\Theta_g$ and $\val(f)=\val(g)$ automatically follow.

\parag{Case (a): $\bs{q = q'}$.}  
(Recall that in this case we have $l=l'$ and $\al_1=\al_1'$.) First we consider the case of $q=q'=1$. In this case by the argument of \cite[Theorem 4.1]{Mi} there exists a repelling cycle $\skakko{\gam_1, \ldots, \gam_{l'}}$ of $f$ satisfying $\gam_{j'} \to \beta_{j'}$ as $f \to g$ and $\Theta(\gam_{j'})=\Theta(\beta_{j'})$ for $j'=1, \ldots, l'$. Take the degenerating arc $\skakko{I_1}$ in the construction of $I_f$. Then $I_1$ joins $\al_1(=\al_1')$ and $\gam_1$ thus $\Theta(\al_1')=\Theta(I(\al_1))=\Theta(\gam_1)=\Theta(\beta_1)$.

Next we consider the case of $q=q'>1$. When $f$ is in segment (s1), the identity $\Theta(\al_1')=\Theta(\beta_1)$ is clear by \lemref{lem_landing}. In the general case, we use renormalization. 

Let us take a path $\eta$ in the parameter space joining $c$ to $\s$ according to the motion as $r \to 1$. By \cite[Th\' eor\` eme 1]{H_Book}, there is an analytic family of quadratic-like maps $\skakko{f^l_{c'}: U_{c'} \to f^l_{c'}(U_{c'})}$ over a neighborhood of $\eta$ such that the straightening maps are continuous and they give one-to-one correspondence between $\eta$ and (s1). 

Let $\al_1 \in O_f$ and $\beta_1 \in O_g$ be the attracting and parabolic fixed points of $f^l=f_c^l:U_{c} \to f^l_{c}(U_{c})$ and $g^l=f_\s^l:U_{\s} \to f^l_{\s}(U_{\s})$ respectively, satisfying $\al_1 \to \beta_1$ as $f \to g$. By \lemref{lem_external}, we can find $q$ external rays landing at $\beta_1$ in the original dynamics of $g$, which is cyclic under $g^l$. In particular, there are no more rays landing at $\beta_1$ since such rays must be cyclic of period $q$ under $g^l$ and this contradicts \cite[Lemma 2.7]{Mi}. Similarly in the dynamics of $f$, by \lemref{lem_landing} and continuity of the straightening, there are exactly $q$ external rays of angles in $\Theta(\beta_1)$ landing at $q$ ends of $I(\al_1)=I(\al_1')$. In fact, if there is another ray of angle $t \notin \Theta(\beta_1)$ landing on such an end, then $R_g(t)$ must land on $\beta_1$ by orbit forcing (\cite[Lemma 7.1]{Mi}). This is a contradiction. Thus $\Theta(\al_1)=\Theta(\al_1')=\Theta(\beta_1)$.

\medskip
\noindent {\bfseries Case (b): $\bs{q = 1<q'}$.}
By the argument of \cite[Theorem 4.1]{Mi}, the repelling points $O'_f=\skakko{\al_1', \ldots, \al'_{l'}}$ must satisfy $\Theta(\al'_{j'})=\Theta(\beta_{j'})$ for $j'=1, \ldots, l'$. 
\QED

\bigskip

In both Cases (a) and (b), it is convenient to assume that $\al_{j'}, \al_{j'+l'}, \ldots, \al_{j'+(q'-1)l'}$ has the same types as that of $\beta_{j'}$ for each $j'=1, \ldots, l'$. Equivalently, we assume that
$$
I(\al_{j'}) \ee  I(\al_{j'+l'}) \ee  \cdots  \ee I(\al_{j'+(q'-1)l'})
$$
throughout this paper. 

\subsection{Critical sectors}\label{subsec_2.3}
For $\xi$ in $I_g$, the external rays of angles in $\Theta(\xi)$ cut the plane up into $\val(g)$ open regions, called \textit{sectors} based at $\xi$. Similarly, for $\zeta$ in $\al_f$, the union of the external rays of angles in $\Theta(\zeta)$ and $I(\zeta)$ cut the plane up into $\val(f)=\val(g)$ open regions. We abuse the term \textit{sectors based at $I(\zeta)$} for these regions. 

Let $B_0$ be the Fatou component of $g$ that contains the critical point $z=0$. We may ssume that $\beta_0=\beta_{l'}$ is on the boundary of $B_0$. Now one of the sectors based at $\beta_0$ contains the critical point $0$, which is called the \textit{critical sector}. For later use, let $\theta_0^+, ~\theta_0^- \in \R/\Z$ denote the angles of external rays bounding the critical sector such that if we take representatives $\theta_0^+ < \theta_0^- \le \theta_0^+ +1$ the external ray of angle $\theta$ with $\theta_0^+ < \theta  <\theta_0^-$ is contained in the critical sector. For example, we define $\theta_0^+:=4/7$ and $\theta_0^-:=1/7$ in the case of \figref{fig_deg_arc_sys}. In the case of \figref{fig_primitive}, we define $\theta_0^+:=5/7$ and $\theta_0^-:=2/7$. We also define the \textit{critical sector} based at $I(\al_0)$ by one of the sector bounded by $I(\al_0)$ and $R_f(\theta_0^\pm)$.


\section{Tessellation}\label{sec_03}
In this section, we develop (and compactify) the method in \cite{Ka3}, and construct \textit{tessellation} of the interior of the filled Julia sets for a degeneration pair $(f \to g)$. 

For each $\theta \in \Theta_f=\Theta_g$ and some $m \in \Z$ (with a condition depending on $\theta$), we will define the tiles $T_f(\theta, m, \pm)$ and $T_g(\theta, m, \pm)$ with the properties listed in \thmref{thm_1}. The idea of tessellation is so simple as one can see in \figref{fig_tessellation}, but we need to construct them precisely to figure out their combinatorial structure in detail.

\subsection{Fundamental model of tessellation}
Take an $R \in (0,1)$ and let us consider the affine maps $F(W)=RW+1$ and $G(W)=W+1$ on $\C$ as the $W$-plane. The map $F$ has a fixed point $a=1/(1-R)$ and one can see the action by the relation $F(W)-a=R(W-a)$.

\parag{Tiles for $\bs{F}$.}
Set $I:=[a, \infty)$, a half line invariant under $F$. For each $\mu \in \Z$, we define ``tiles" of level $\mu$ for $F$ by:
\begin{align*} 
A_\mu (+)&:=
\skakko{W \in \C-I: R^{\mu +1}a \le |W-a| \le R^\mu a
,~\ip{W} \ge 0}\\
A_\mu (-)&:=
\skakko{W \in \C-I: R^{\mu +1}a \le |W-a| \le R^\mu a
,~\ip{W} \le 0}.
\end{align*}
Then one can check that $F(A_\mu (\ast))=A_{\mu +1}(\ast)$ where $\ast \in \skakko{+, -}$. For the boundary of each $A_\mu (\ast)$, we define 
\begin{itemize}
\item the \textit{circular edges} by the intersection with $A_{\mu  \pm 1}(\ast)$;
\item the \textit{degenerating edge} by $\overline{A_\mu (\ast)} \cap I$; and
\item the \textit{critical edge} by the intersection with $(-\infty, a)$. 
\end{itemize}
Note that $A_\mu (\ast) \subset \C-I$ so the degenerating edge is \textit{not} contained in $A_\mu (\ast)$.

\parag{Tiles for $\bs{G}$.}
Correspondingly, for each $\mu  \in \Z$, we define ``tiles" of level $\mu$ for $G$ by:
\begin{align*} 
C_\mu (+)&\dee \skakko{W \in \C: \mu  \le \rp{W} \le \mu +1
,~\ip{W} \ge 0}\\
C_\mu (-)&\dee\skakko{W \in \C: \mu  \le \rp{W} \le \mu +1
,~\ip{W} \le 0}.  
\end{align*}
Then one can check that $G(C_\mu (\ast))=C_{\mu +1}(\ast)$. For the boundary of each $C_\mu (\ast)$, we define 
\begin{itemize}
\item the \textit{circular edges} by the intersection with $C_{\mu \pm 1}(\ast)$, which are vertical half lines;
\item the \textit{critical edge} by the intersection with $(-\infty, \infty)$. 
\end{itemize}
Note that there is no degenerating edge for $C_\mu (\ast)$. One can consider $C_\mu (\ast)$ the limit of $A_\mu (\ast)$ as $R \to 1$.

\begin{figure}[htbp]
\centering
\includegraphics[width=0.85\textwidth]{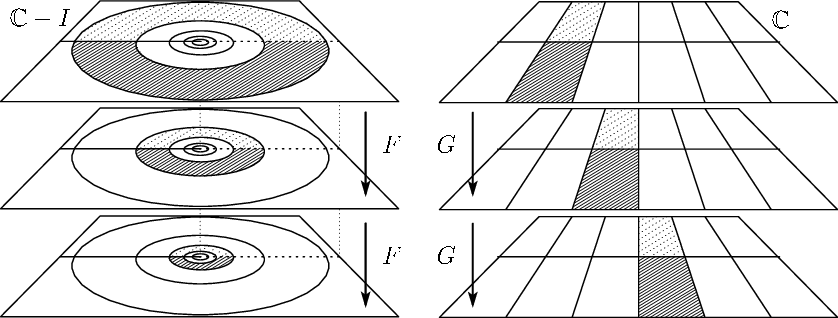}
\caption{The fundamental model}\label{fig_FG}
\end{figure}

\subsection{Tessellation for $\bs{f}$ and $\bs{g}$.}

First we reduce the dynamics of $f|_{\Kfc}$ and $g|_{\Kgc}$ to the dynamics of $F$ and $G$ on $\C$.

\parag{From $\bs{f}$ to $\bs{F}$.}
Let $B_0$ be the Fatou component of $f$ containing $0$. We may assume that $\al_0=\al_l \in B_0$. There exists a unique extended linearizing coordinate $\phi_f:B_0 \to \C$ such that $\phi_f(\al_0)=\phi_f(0)-1=0$ and $\phi_f(f^{l}(z)) = \lam \phi_f(z)$ \cite[\S 8]{MiBook}. Set $w:=\phi_f(z)$ and $R:=\lambda^q =r^q$. Then $f^{lq}|_{B_0}$ is semiconjugate to $w \mapsto Rw$. To reduce this situation to our fundamental model, first we take a branched covering $W=w^q$. Then $f^l|_{B_0}$ and $f^{lq}|_{B_0}$ are semiconjugate to $W \mapsto RW$ and $W \mapsto R^qW$ respectively. Next, we take an affine conjugation by $W \mapsto a(1-W)$. Then $f^l|_{B_0}$ and $f^{lq}|_{B_0}$ are finally semiconjugate to $F$ and $F^q$ in the fundamental model respectively. Let $\Phi_f$ denote this final semiconjugation. Now we have $\Phi_f(0)=0$ and $\Phi_f(B_0 \cap I_f)=I$. (The second equality comes from the construction of the degenerating arcs in \lemref{lem_landing}.) In particular, the map $\Phi_f$ branches at $z \in B_0$ if and only if either $f^{ln}(z)=0$ for some $n \ge 0$, or $q>1$ and $f^{ln}(z)=\al_0$ for some $n \ge 0$.

\parag{From $\bs{g}$ to $\bs{G}$.}
Let $B_0'$ be the Fatou component of $g$ containing $0$. We may assume that $\beta_0=\beta_{l'} \in \partial B_0'$. There exists a unique extended Fatou coordinate $\phi_g:B_0' \to \C$ such that $\phi_g(0)=0$ and $\phi_g(g^{lq}(z)) = \phi_g(z)+1$ \cite[\S 10]{MiBook}. Set $w:=\phi_g(z)$, then $g^{lq}|_{B_0}$ is semiconjugate to $w \mapsto w+1$. To adjust the situation to that of $f$, we take an additional conjugacy by $w \mapsto W=qw$. Then $g^{lq}|_{B_0}$ is semiconjugate to $G^q(W)=W+q$. We denote this semiconjugation $z \mapsto w \mapsto W$ by $\Phi_g$. Note that $\Phi_g(0)=0$, and $\Phi_g$ branches at $z \in B_0'$ iff $g^{lqn}(z)=0$ for some $n \ge 0$.

\begin{figure}[htbp]
\vspace{0cm}\hspace{0cm}
\centering{
\includegraphics[height=0.4\textheight]{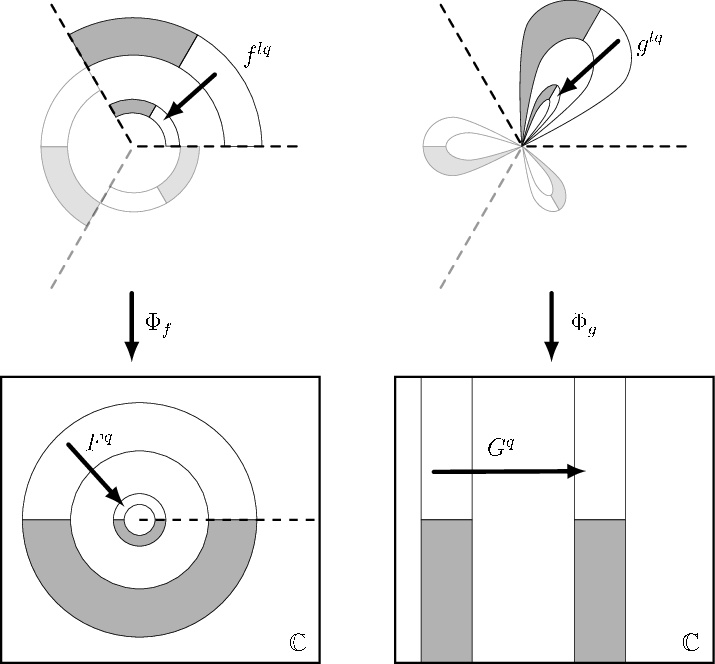}
}
\caption{$f^{lq}$ and $g^{lq}$ are semiconjugate to $F^q$ and $G^q$.}\label{fig_adjusting_levels}
\end{figure}

\medbreak
Let us summarize these reduction steps. Now $\Phi_f:B_0 \to \C$ semiconjugates the action of $f^{lq}: B_0-I_f \to B_0-I_f$ to that of $F^q:\C-I \to \C -I$. Similarly, the map $\Phi_g:B_0' \to \C$ semiconjugates the action of $g^{lq}: B_0' \to B_0'$ to that of $G^q:\C \to \C$ (\figref{fig_adjusting_levels}). In addition, we have one important property as follows:

\begin{prop}\label{prop_ramify}
The branched linearization $\Phi_f$ do not ramify over $\C-(-\infty, 0]$ or $\C-(-\infty,0]\cup\skakko{a}$ according to $q=1$ or $q>1$. Similarly, the branched linearization $\Phi_g$ do not ramify over $\C-(-\infty, 0]$. In particular, both $\Phi_f$ and $\Phi_g$ do not ramify over tiles of level $\mu>0$.
\end{prop}
See \thmref{thm_linearization_converges} for another important property of $\Phi_f$ and $\Phi_g$.

\parag{Definition of tiles.}
A subset $T \subset \Kfc$ is a \textit{tile} for $f$ if there exist $n \in \N$ and $\mu  \in \Z$ such that $f^n(T)$ is contained in $B_0$ and $\Phi_f \cc f^n$ maps $T$ homeomorphically onto $A_\mu (+)$ or $A_\mu (-)$. We define \textit{circular, degenerating, and critical edges} for $T$ by their corresponding edges of $A_\mu (\pm)$. We call the collection of such tiles the \textit{tessellation} of $\Kfc-I_f$, and denote it by $\Tess(f)$. In fact, one can easily check that
$$
\Kfc-I_f \ee \bigcup_{T \in \Tess(f)} T
$$
and each $z \in \Kfc -I_f$ is either in the interior of an unique $T \in \Tess(f)$; a vertex shared by four or eight tiles in $\Tess(f)$ if $f^m(z)=f^n(0)$ for some $n,m >0$; or on an edge shared by two tiles in $\Tess(f)$ otherwise.

Tiles for $g$ and tessellation of $\Kgc-I_g=\Kgc$ are also defined by replacing $f$, $B_0$, and $A_\mu (\pm)$ by $g$, $B_0'$, and $C_\mu (\pm)$ respectively.

\parag{Addresses.} 
Each tile is identified by an \textit{address}, which consists of an \textit{angle}, a \textit{level}, and a \textit{signature} defined as followings:

\parag{Level and signature.} 
For $T \in \Tess(f)$ above, i.e., $f^n(T) \subset B_0$ and $\Phi_f \cc f^n(T)=A_\mu (\ast)$ with $\ast=+$ or $-$, we say that $T$ has \textit{level} $m=\mu l-n$ and \textit{signature} $\ast$. Then the critical point $z=0$ is a vertex of eight tiles of level $0$ and $-l$.

For a tile $T' \in \Tess(g)$, its level and signature is defined in the same way.

\begin{figure}[htbp]
\vspace{0cm}\hspace{0cm}
\centering{
\includegraphics[width=0.9\textwidth]{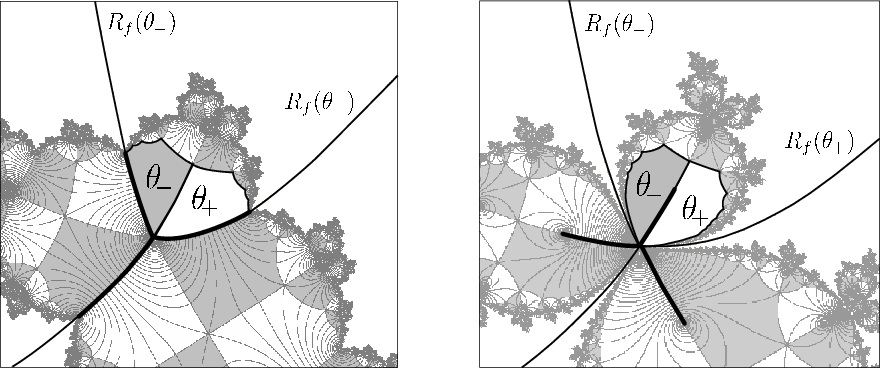}
}
\caption{Angles of tiles in Case (a) (left) and Case (b) (right) with $q'=3$. The thick arcs show degenerating arcs.}
\label{fig_angle}
\end{figure}

\parag{Angles.}
For $T \in \Tess(f)$, there exists $\zeta$ in $\al_f$ such that $I(\zeta)$ contains the degenerating edge of $T$. Then there are $\val(f)=v \ge 1$ rays landing on $I(\zeta)$, and the rays and $I(\zeta)$ divide the plane into $v$ sectors. (In the case of $v=1$, equivalently $g(z)=z^2+1/4$, we consider the sector as the plane with a slit.) Take two angles $\theta_{+} < \theta_{-} (\leq \theta_{+} + 1)$ of external rays bounding the sector containing $T$. (That is, any external ray of angle $\theta$ with $\theta_{+} < \theta < \theta_{-}$ is contained in the sector.) Now we define the \textit{angle} of $T$ by $\theta_{\ast}$ where $\ast$ is the signature of $T$. (See \figref{fig_angle})

For a tile $T' \in \Tess(g)$, one can check that there exists a unique point $\beta' \in I_g \cap \partial T'$. Since there are $v$ rays land on $\beta'$ and they divide the plane into $v$ sectors as in the case of $T \in \Tess(f)$, we define the angle of $T'$ in the same way as above.

\begin{figure}[htbp]
\centering
\vspace{-1cm}\hspace{0cm}
\includegraphics[height=0.9\textheight]{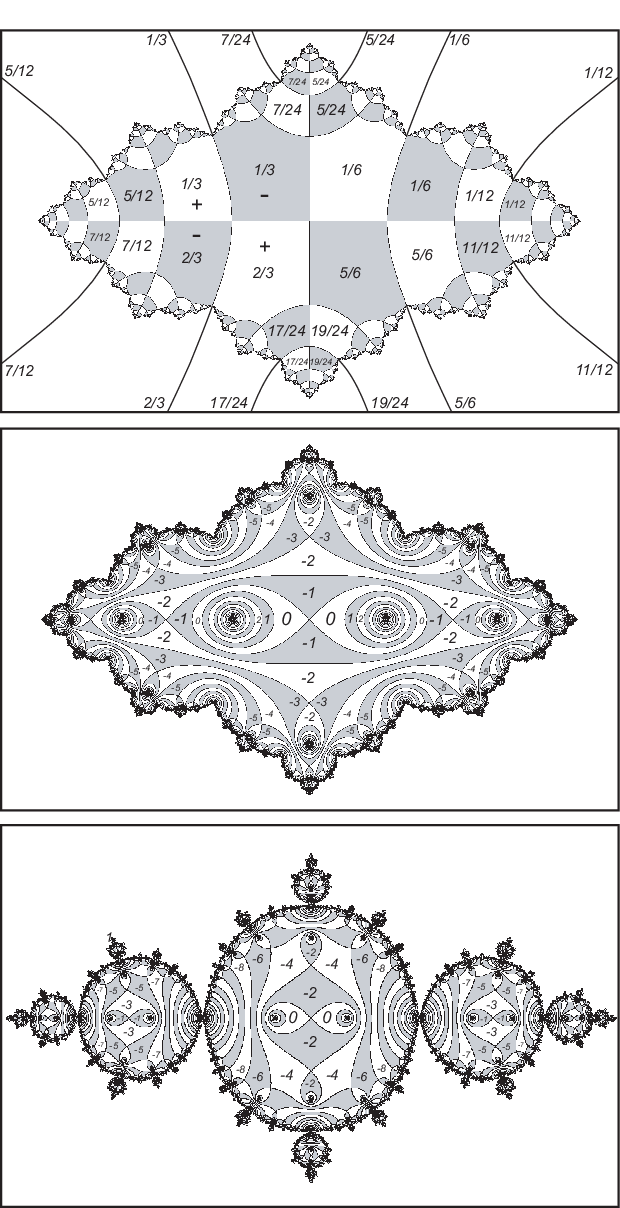}
\caption{``Checkerboard" and ``Zebras", showing the structure of the addresses of tiles. ``Checkerboard", with some external rays drawn in, shows the relation between the external angles and the angles of tiles. The invariant regions colored in white and gray correspond to tiles of signature $+$ and $-$ respectively. ``Zebras" show the levels of tiles for $f_c$ with $-1 < c < 0$. Levels get higher near the preimages of the attracting periodic points.}\label{fig_check_zebra}
\end{figure}

We denote such tiles by $T=T_f(\theta_\ast,m,\ast)$ and $T'=T_g(\theta_\ast,m,\ast)$, and we call the triple $(\theta_\ast, m, \ast)$ the \textit{address} of the tiles. For example, \figref{fig_check_zebra} shows the structure of addresses for the two tessellations at the lower left of \figref{fig_tessellation}. 

Now one can easily check the desired property
$$
f(T_f(\theta, m, \ast)) \ee T_f(2\theta, m+1, \ast).
$$
The same holds if we replace $f$ by $g$. One can also check properties (1) to (5) of \thmref{thm_1} easily.

\begin{figure}[htbp]
\centering
\hspace*{0cm}\vspace*{0cm}
\includegraphics[width=0.95\textwidth]{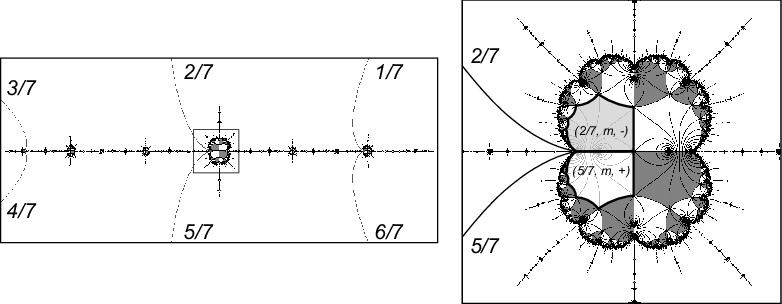}
\caption{A non-trivial primitive $(f \to g)$ with $g(z)=z^2-7/4$. For example, we define $T_f(2/7,m,+)$ by $T_f(5/7,m,+)$.}
\label{fig_primitive}
\end{figure}

\parag{Remarks on angles and levels.}
\begin{itemize}
\item We make an exception for non-trivial primitives ($q=q'$ and $l=l'>1$). If $(f \to g)$ is non-trivial primitive, then $v=2$ and only tiles of addresses $(\theta_\pm, m, \pm)$ are defined. However, we formally define tiles of addresses $(\theta_\pm, m, \mp)$ by tiles of addresses $(\theta_\mp, m,  \mp)$ respectively. (See \figref{fig_primitive}.)

\item For a degeneration pair $(f \to g)$, the space of possible addresses of tiles is not equal to $\Theta_f \times \Z \times \skakko{+, -}$ in general. For both $f$ and $g$, all possible addresses are realized when $l=1$. But when $l > 1$, the address $(\theta, m, \pm)$ is realized iff $m+n \equiv 0 \mod l$ for some $n>0$ with $2^n \theta =\theta_0^\pm$. \textit{In any case, note that $T_f(\theta, m, *)$ exists iff $T_g(\theta, m, *)$ exists.}
\end{itemize}

\subsection{Edge sharing}
Let us investigate the combinatorics of tiles in $\Tess(f)$ and $\Tess(g)$. We will show the following proposition that is a detailed version of \thmref{thm_1}(6):

\begin{prop}\label{prop_edge_sharing}
For $\theta \in \Theta_f=\Theta_g$ and $\ast \in \skakko{+, -}$, let us take an $m \in \Z$ such that $T=T_f(\theta, m , \ast)$ and $S=T_g(\theta, m, \ast)$ exist. Then:
\begin{enumerate}
\item The circular edges of $T$ and $S$ are shared by $T_f(\theta, m \pm l, \ast)$ and $T_g(\theta, m \pm l, \ast)$ respectively. 
\item The degenerating edge of $T$ is contained in $I(\zeta)$ with $\zeta \in \al_f$ of type $\Theta(\zeta)$ iff $S$ attaches at $\xi \in I_g$ of type $\Theta(\xi)=\Theta(\zeta)$. Moreover, the degenerating edge of $T$ is shared with $T_f(\theta, m, \bar{\ast})$, where $\bar{\ast}$ is the opposite signature of $\ast$.
\item $T$ shares its critical edge with $T_f(\theta', m', \ast')$ iff $S$ does the same with $T_g(\theta', m', \ast')$. In this case, we have $m'=m$ and $\ast'=\bar{\ast}$. 
\end{enumerate}
\end{prop}
Thus the combinatorics of $\Tess(f)$ and $\Tess(g)$ are the same.

\parag{Proof. (1) Circular edges:}
By \propref{prop_ramify}, for any $n \ge 0$, the inverse map $f^{-n} \cc \Phi_f^{-1}$ over $\C-(-\infty, a]$ is a multi-valued function with univalent branches. Now it follows that the property ``$A_\mu(\ast)$ shares its circular edges with $A_{\mu \pm 1}(\ast)$" is translated to ``$T(\theta, m, \ast)$ shares its circular edges with $T(\theta, m \pm l, \ast)$" by one of such univalent branches. The same argument works for $\Phi_g:B_0' \to \C$, which does not ramify over $\C-(-\infty,0]$.

\parag{(2) Degenerating edges:} The statement is clear by definition of tiles and addresses.

\parag{(3) Critical edges:} 
The combinatorics of tiles are essentially determined by the connection of critical edges. They are organized as follows.

In the fundamental model, we consider a family of curves
\begin{align*}
\overline{A_\mu (\ast)} \cap \skakko{|W-a|=R^{\mu +1/2}} \\
C_\mu (\ast) \cap \skakko{\rp W=\mu R+1/2}
\end{align*}
for $\mu  \in \Z$ and we call the \textit{essential curves} of $A_\mu (\ast)$ and $C_\mu (\ast)$. Since $\Phi_f \cc f^n$ and $\Phi_g \cc g^n$ do not ramify over these essential curves, their pulled-back images in the original dynamics form ``equipotential curves" in $\Kfc$ and $\Kgc$. The \textit{essential curve} of a tile is the intersection with such equipotential curves.

Let us consider a general tile $T \in \Tess(f)$ as in the statement. By taking a suitable $n \gg 0$, we may assume that $f^n(T)$ is a tile in $B_0$ with angle $t$ in $\skakko{\theta_0^+, \theta_0^-} \subset \Theta(\al_0)$ and level $\mu l$ for some $\mu \ge 0$. In particular, we may assume that $f^n(T)$ is in the critical sector based at $I(\al_0)$. Then for $S$ in the statement, we can take the same $n$ and $\mu$ as those for $T$ such that $g^n(S)$ is a tile in $B_0'$ with angle $t$ in $\Theta(\beta_0)=\Theta(\al_0)$ and level $\mu l$.

\parag{Case (a): $\bs{q =q'}$.} 
Let $\eta_0$ be the union of essential curves of tiles of the form $T_f(t, \mu l, \ast)$ with $t$ in $\Theta(\al_0)$. Then $\eta_0$ forms an equipotential curve around $\al_0$, since $\Phi_f|_{\eta_0}$ is a $q$-fold covering over the circle $\skakko{|W-a|=R^{1/2+\mu}}$. 
For $n > 0$, set $\eta_{-n}=f^{-n}(\eta_0)$. Then $\eta_{-n}$ is a disjoint union of simple closed curves passing through tiles of level $\mu l - n$ and angles in $\delta^{-n}(\Theta(\al_0))$. In particular, each curve crosses degenerating edges and critical edges alternatively. More precisely, let $\eta$ be a connected component of $\eta_{-n}$. Then the degree of $f^n:\eta \to \eta_0$ varies according to how many curves in $\skakko{f^k(\eta)}_{k=1}^n$ enclose the critical point $z=0$. One can check the degree by counting the number of points of $f^{-n}(\al_0)$ inside $\eta$. Let $\zeta_1, \ldots, \zeta_N$ be such points. Then $\eta$ crosses each $I(\zeta_i)$, and thus $\eta$ crosses the tiles of level $-n$ with angles in $\Theta(\zeta_1) \cup \cdots \cup \Theta(\zeta_N) \subset \T$ in cyclic order, and with signatures switching as crossing the edges of tiles. This observation gives us how critical and degenerating edges are shared among tiles along $\eta$. 

Now we can take $\eta$ passing through $T$. This observation concludes: if $T$ shares its critical edge with $T_f(\theta', m', \ast')$, then $m'=m$, and $\ast'=\bar{\ast}$; and if $T$ shares its degenerating edge with $T_f(\theta', m', \ast')$, then $\theta'=\theta$, $m'=m$, and $\ast'=\bar{\ast}$.

\begin{figure}[htbp]
\vspace{0cm}\hspace{0cm}
\centering{
\includegraphics[width=0.7\textwidth]{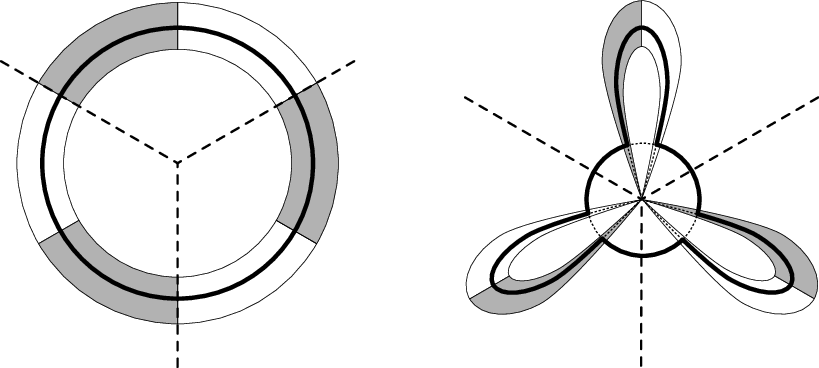}
}
\caption{The thick curves show $\eta_0$ and $\eta_0'$ in Case (a) with $q=q'=3$. The dashed lines indicate the degenerating arcs or external rays.}\label{fig_eta_0_case_1}
\end{figure}

For $S$, consider a circle around $\beta_0$ which is so small that the circle and the essential curves of tiles with angle $\theta \in \Theta(\beta_0)$ and level $\mu l$ bound a flower-like disk (\figref{fig_eta_0_case_1}). Let us denote the boundary of the disk by $\eta'_0$, which works as $\eta_0$. Since the combinatorics of pulled-back sectors based at $\beta_0$ and $I(\al_0)$ is the same, the observation of $g^{-n}(\eta_0')=\eta'_{-n}$ must be the same as that of $\eta_{-n}$. This concludes the statement.

\parag{Case (b): $\bs{q = 1<q'}$.} 
(Recall that in this case $O_g$ is perturbed into the repelling cycle $O_f'=\skakko{\al'_{1}, \ldots, \al'_{l'}=\al'_0}$ with $\al'_0 \to \beta_0$ as $f \to g$.)
The same argument as above works if we take $\eta_0$ and $\eta_0'$ as following: First in the fundamental model, take $\e \ll 1$ and two radial half-lines from $a$ with argument $\pm \e$. Then there are univalently pulled-back arcs of two lines in the critical sector which joins $\al_0$ and $\al'_0$. Next we take simple closed curves around $\al_0$ and $\al'_0$. For $\al_0$, we take the essential curves along tiles of address $(\theta_0^\pm, \mu l, \pm)$. For $\al'_0$, we take just a small circle around $\al'_0$. Then the two arcs and two simple closed curves bound a dumbbell-like topological disk. We define $\eta_0$ by its boundary curve.

\begin{figure}[htbp]
\vspace{0cm}\hspace{0cm}
\centering{
\includegraphics[width=0.7\textwidth]{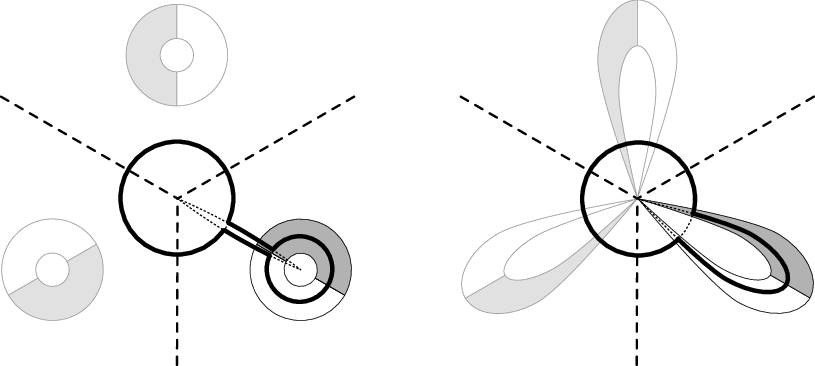}
}
\caption{$\eta_0$ and $\eta_0'$ in Case (b) with $q'=3$.}
\label{fig_eta_0_case_2}
\end{figure}

Correspondingly, for $g$, the essential curves of tiles of address $(\theta_0^\pm, \mu l, \pm)$ and a small circle around $\beta_0$ bound a topological disk. We take $\eta'_0$ as its boundary (\figref{fig_eta_0_case_2}).
\QED

\subsection{Tiles and panels with small diameters}
Next we show that the diameter of tiles are controlled by their angles. For $\theta$ in $\Theta_f=\Theta_g$ and $\ast \in \skakko{+,-}$, let $\Pi_f(\theta, \ast)$ and $\Pi_g(\theta, \ast)$ be the union of tiles with angle $\theta$ and signature $\ast$ in $\Tess(f)$ and $\Tess(g)$ respectively. We call them \textit{panels} of angle $\theta$ and signature $\ast$. (For later use, by $\Pi_f(\theta)$ we denote $\Pi_f(\theta, +) \cup \Pi_f(\theta, -)$.) The \textit{depth} of angle $\theta$ is the minimal $n \ge 0$ such that 
$
2^n \theta = \theta_0^+,
$
where $\theta_0^+ \in \Theta(\al_0)=\Theta(\beta_0)$ is defined in Subsection \ref{subsec_2.3}. (Note that $\Pi_f(\theta_0^+)=\Pi_f(\theta_0^-)$ when $(f \to g)$ is non-trivial primitive.) We denote such an $n$ by $\depth(\theta)$. Here we show the following:

\begin{prop}\label{prop_diam}
For any fixed degeneration pair $(f \to g)$ and any $\e >0$, there exists $N=N(\e, f, g)$ such that 
$$
\diam \Pi_f(\theta, \ast) < \e ~~~\text{and}~~~\diam \Pi_g(\theta, \ast) < \e 
$$
for any signature $\ast$ and any $\theta \in \Theta_f$ with $\depth(\theta) \ge N$.
\end{prop}

\parag{Proof.}
We first work with $f$ and signature $+$. One can easily check that the interior $\Pi$ of $\Pi_f(\theta_0^+,+)$ is a topological disk. For any $\theta \in \Theta_f$, the panel $\Pi_f(\theta,+)^\cc$ is sent univalently onto $\Pi$ by $f^n$ with $n=\depth(\theta)$. Let $F_\theta$ be the univalent branch of $f^{-n}$ which sends $\Pi$ to $\Pi_f(\theta,+)^\cc$. Since the family $\skakko{F_\theta:\theta \in \Theta_f}$ on $\Pi$ avoids the values outside the Julia set, it is normal. 

Now we claim: \textit{$\diam \Pi_f(\theta,+) \to 0$ as $\depth(\theta) \to \infty$}. Otherwise one can find a sequence $\skakko{\theta_k}_{k >0}$ with depth $n_k \to \infty$ and $\delta >0$ such that $\diam \Pi_f(\theta_k,+) >\delta$ for any $k$. By passing through a subsequence, we may assume that $F_k:=F_{\theta_k}$ has a non-constant limit $\phi$. Fix any point $z \in \Pi$, and set $\zeta:=\phi(z)=\lim F_k(z)$. Since $\phi$ is holomorphic and thus is an open map, there exists a neighborhood $V$ of $\zeta$ such that $V \subset \phi(\Pi)$ and $V \subset F_k(\Pi)$ for all $k \gg 0$. Since $f^{n_k}(V) \subset \Pi \subset \Kfc$, any point in $V$ are attracted to the cycle $O_f$. However, by univalence of $F_k$, there exists a neighborhood $W$ of $z$ with $W \subset F_k^{-1}(V)=f^{n_k}(V)$ for all $k \gg 0$. This is a contradiction. 

Finally we arrange the angles of $\Theta_f$ in a sequence $\skakko{\theta_i}_{i >0}$ such that $\depth(\theta_n)$ is non-decreasing. Note that for any integer $n$, there are only finitely many angles with depth $n$. Thus there exists an integer $N=N(\e, f,+)$ such that $\Pi_f(\theta, +)$ has diameter less than $\e$ if $\depth(\theta) \ge N$. 

This argument works if we switch the map (from $f$ to $g$) or the signature. Then we have four distinct $N$ as above. Now we can take $N(\e, f, g)$ as their maximum.
\QED

\paragraph{}
Indeed, as depth tends to infinity we have uniformly small panels for $f \approx g$ (\propref{prop_uniformly_small_panels}).

\section{Pinching semiconjugacy}\label{sec_04}
In this section we construct a semiconjugacy $h: \Cbar \to \Cbar$ associated with $(f \to g)$ by gluing tile-to-tile homeomorphisms inside the Julia sets and the topological conjugacy induced from the B\" ottcher coordinates outside the Julia sets.

\begin{thm} \label{thm_semiconj} 
For a degeneration pair $(f \to g)$, there exists a semiconjugacy $h:\Cbar \to \Cbar$ from $f$ to $g$ such that
\begin{enumerate}[\rm (1)] 
\item $h$ maps $\Cbar-I_f$ to $\Cbar-I_g$ homeomorphically and is a topological conjugacy between $f|_{\Cbar-I_f}$ and $g|_{\Cbar-I_g}$;  
\item For each $\zeta \in \al_f$ with type $\Theta(\zeta)$, $h$ maps $I(\zeta)$ onto a point $\xi \in I_g$ with type $\Theta(\xi)=\Theta(\zeta)$.
\item $h$ sends all possible $T_f(\theta, m, \ast)$ to $T_g(\theta, m, \ast)$, $R_f(\theta)$ to $R_g(\theta)$, and $\gam_f(\theta)$ to $\gam_g(\theta)$.
\end{enumerate}
\end{thm}

This theorem emphasizes the combinatorial property of $h$. In the next section we will show that $h \to \id$ as $f$ uniformly tends to $g$. 

\parag{Trans-component partial conjugacy and subdivision of tessellation.} 
Let $(f_1 \to g_1)$ and $(f_2 \to g_2)$ be distinct satellite degeneration pair with $g_1=g_2$. More precisely, we consider $(f_1 \to g_1)$ and $(f_2 \to g_2)$ are tuned copy of degeneration pairs in segment (s1) and (s2) with $q>1$ by the same tuning operator. By composing homeomorphic parts of the conjugacies associated with $(f_1 \to g_1)$ and $(f_2 \to g_2)$, we have:

\begin{cor} \label{cor_trans_compo} 
There exists a topological conjugacy $\kappa=\kappa_{f_1,f_2}:\Cbar-I_{f_1} \to \Cbar-I_{f_2}$ from $f_1$ to $f_2$.
\end{cor}

For example, the panel $\Pi_{f_1}(\theta, \ast)$ is mapped to the panel $\Pi_{f_2}(\theta, \ast)$. Now we can compare $\Tess(f_1)$ and $\Tess(f_2)$ via $\Tess(g_i)$. By comparing $\Tess(g_1)$ and $\Tess(g_2)$, one can easily check that 
$$
T_{g_2}(\theta, \mu, \ast)
 \ee \bigcup_{j=0}^{q-1}T_{g_1}(\theta,~ \mu + lj,~ \ast)
$$
for any $T_{g_2}(\theta, \mu, \ast) \in \Tess(g_2)$. Thus $\Tess(g_1)$ is just a subdivision of $\Tess(g_2)$. 

Take a tile $T_{f_1}(\theta, m , \ast) \in \Tess(f_1)$. Then there is a homeomorphic image $T_{2}'(\theta, m , \ast):=\kappa(T_{f_1}(\theta, m , \ast))$ in $K_{f_2}^\cc$. We say the family
$$
\Tess'(f_2):=\skakko{\kappa(T): T \in \Tess(f_1)}
$$
is the \textit{subdivided tessellation} of $K_{f_2}^\cc-I_{f_2}$. Since $\Tess(f_1)$ and $\Tess(f_2)$ have the same combinatorics as $\Tess(g_1)$ and $\Tess(g_2)$ respectively,
$$
T_{f_2}(\theta, \mu, \ast)
\ee \bigcup_{j=0}^{q-1}T_{f_2}'(\theta, ~\mu+lj,~ \ast)
$$
for any $T_{f_2}(\theta, \mu, \ast) \in \Tess(f_2)$. Now we have a natural tile-to-tile correspondence between $\Tess(f_1)$, $\Tess(g_1)$ and $\Tess'(f_2)$. In other word, combinatorial property of tessellation is preserved under the degeneration from $f_1$ to $g$ and the bifurcation from $g$ to $f_2$.

In Part II of this paper, we will use this property to investigate the structures of the Lyubich-Minsky hyperbolic 3-laminations associated with $f_1$, $g$, and $f_2$.

\subsection{Proof of \thmref{thm_semiconj}}
The rest of this section is devoted to the proof of the theorem. The proof breaks into five steps.

\parag{1. Conjugacy on the fundamental model.}
First we define a topological map $H:\C-I \to \C$ which maps $A_\mu(\pm)$ to $C_\mu(\pm)$ homeomorphically. For $W \in \C-I$, set $W:=a+\rho e^{it}$ where $\rho>0$ and $0 <t< 2 \pi$. We define the map $H$ by
$$
H(W):=\frac{\log \rho-\log a}{\log R}+ 2a i \tan\frac{\pi-t}{2} \in \C.
$$
Then one can check that $H$ conjugates the action of $F$ on $\C-I$ to that of $G$ on $\C$ and $H$ maps $A_\mu(\pm)$ homeomorphically onto $C_\mu(\pm)$. 

\begin{figure}[htbp]
\centering
\includegraphics[width=0.85\textwidth]{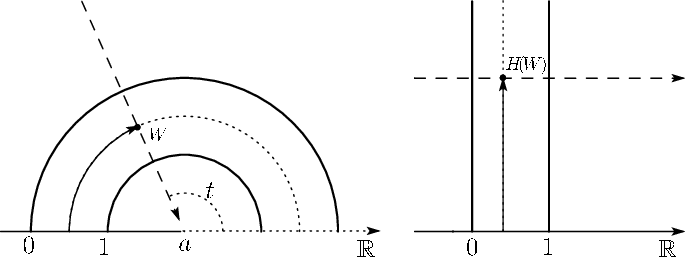}
\caption{$H$ maps $A_0(+)$ to $C_0(+)$.}
\end{figure}

\parag{2. Tile-to-tile conjugation.}
First we consider the critical sectors of $f$ and $g$. Let $\Pi_0$ and $\Pi_0'$ denote the union of tiles of addresses $(\theta_0^\pm, \mu l, \pm)$ with $\mu > 0 $ in $\Tess(f)$ and $\Tess(g)$ respectively.

By \propref{prop_ramify}, the map $\Phi_f:\Pi_0 \to \C$ is univalent and we can choose a univalent branch $\Psi_g$ of $\Phi_g^{-1}$ which sends $\skakko{W: \rp W \ge 1}$ to $\Pi_0'$. For each point in $\Pi_0$, we define $h:=\Psi_g \cc H \cc \Phi_f|_{\Pi_0}$. Then $h$ is a conjugacy between $f^{lq}|_{\Pi_0}$ and $g^{lq}|_{\Pi_0'}$. Note that any tiles eventually land on tiles in $\Pi_0$ or $\Pi'_0$. According to the combinatorics of tiles determined by pulling-back essential curves in $\Pi_0$ and $\Pi_0'$, we can pull-back $h$ over $\Kfc-I_f$. Now the map $h:\Kfc-I_f \to \Kgc$ conjugates $f|_{\Kfc-I_f}$ and $g|_{\Kgc}$.

\parag{3. Continuous extension to degenerating arc system.}
Take any $\zeta \in \al_f$.  For any point $z$ in $I(\zeta)$, we define $h(z)$ by the unique $\xi \in I_g$ with $\Theta(\xi)=\Theta(\zeta)$. 

Now we show the continuity of $h:\Kfc \cup I_f \to \Kgc \cup I_g$ which we have defined. Take any $z$ in $I(\zeta)$. We claim that any sequence $z_n \in \Kfc \cup I_f$ converging to $z$ satisfies $h(z_n) \to \xi$. 

First when $z$ is neither $\zeta$ nor one of the endpoints of $I(\zeta)$, it is enough to consider the case of $z_n \in \Kfc - I_f$ for all $n$. Now $z$ is on the degenerating edges of at most four tiles. Let $T=T_f(\theta, m, +)$ be one of such tiles. The subsequence $z_{n_i}$ of $z_n$ contained in $T$ is mapped to $T_g(\theta, m, +)$. In the fundamental model, the sequence $h(z_{n_i})$ corresponds to a sequence whose imaginary part is getting higher. Thus $h(z_{n_i})$ converges to $\xi$ with type containing $\theta$, which must coincide with $\Theta(\zeta)$. By changing the choice of $T$, we have $h(z_n) \to \xi$ with $\Theta(\xi)=\Theta(\zeta)$.

Next, if $z$ is $\zeta$ or one of the endpoints of $I(\zeta)$, it is an attracting or repelling periodic point. If $z$ is attracting, the levels of tiles containing $z_n$ go to $+ \infty$. According to the fundamental model, we have $h(z_n) \to \xi$.

The remaining case is that $z$ is repelling, thus in the Julia set. We deal with this case in the next paragraph.

\parag{4. Continuous extension to the Julia set.}
Take any $z \in J_f$ and any sequence $z_n \in \Kfc \cup I_f$ converging to $z$. Then we take a sequence $\theta_n \in \Theta_f$ such that $z_n \in \Pi_f(\theta_n)$. After passing to a subsequence we may assume $\theta_n$ and $h(z_n)$ converge to some $\theta \in \T$ and $w \in K_g$ respectively.

We first claim that $z=\gam_f(\theta)$, that is, $\theta \in \Theta(z)$. If the depth of $\theta_n$ is bounded, then $\theta_n = \theta \in \Theta_f$ for all $n \gg 0$. 
This implies $z_n \in \Pi_f(\theta)$ for all $n \gg 0$ and it follows that $z \in \overline{\Pi_f(\theta)} \cap J_f$. Thus $z =\gam_f(\theta)$ by definition of $\Pi_f(\theta)$. 
If the depth of $\theta_n$ is unbounded, it is enough to consider the subsequences with the depth of $\theta_n$ monotonously increasing. 
Take any $\e >0$. 
For $n \gg 0$, we have $|\gam_f(\theta_n)-z_n|<\e$ by \propref{prop_diam}, and we also have $|\gam_f(\theta_n)-\gam_f(\theta)|<\e$ by continuity of $\gam_f:\T \to J_f$. 
Finally $|z - z_n| < \e$ for $n \gg 0$ implies $|z - \gam_f(\theta)| < 3 \e$ and we conclude the claim.

Since $h(z_n) \in \Pi_g(\theta_n)$, the same argument works for $h(z_n)$ and $w$. Hence we also have $w = \gam_g(\theta) \in J_g$. It follows that for the original $z_n \to z$, the sequence $h(z_n)$ accumulates only on $\gam_g(\theta)$ with $\theta \in \Theta(z)$. 

By \thmref{thm_semiconj_on_julia}, there exists a semiconjugacy $h_J:J_f \to J_g$ with $h_J \cc \gam_f = \gam_g$. Since $\gam_f(\theta)=\gam_f(\theta')$ for any $\theta$ and $\theta'$ in $\Theta(z)$, we have $\gam_g(\theta)=\gam_g(\theta')$. This implies that $h(z_n)$ accumulates on a unique point $\gam_g(\theta)$. Thus $h$ continuously extends to the Julia set by $h(\gam_f(\theta)):=\gam_g(\theta)$ for each $\theta \in \T$.





\parag{5. Global extension.}
Finally we define $h:\Cbar-K_f \to \Cbar -K_g$ by the conformal conjugacy between $f|_{\Cbar-K_f}$ and $g|_{\Cbar-K_g}$ given via the B\" ottcher coordinates. This conjugacy and the semiconjugacy above continuously glued along the Julia set thus we have a semiconjugacy on the sphere.

Properties (2) and (3) are clear by construction. To check property (1), we need to show that $h^{-1}: \Cbar-I_f \to \Cbar - I_g$ is continuous. Continuity in $\Cbar-K_g$ and $\Kgc$ is obvious by construction. Take any point $w \in J_g-I_g$. A similar argument to step 4 shows that any sequence $w_n \to w$ within $\Cbar-I_g$ is mapped to a convergent sequence $z_n \to z$ within $\Cbar-I_f$ satisfying $\Theta(z) = \Theta(w) \subset \T-\Theta_g$. 
\QED

\section{Continuity of pinching semiconjugacies}\label{sec_05}
In this section we deal with continuity of the dynamics of the degeneration pair $(f \to g)$  as $f$ tends to $g$. We will establish:

\begin{thm}\label{thm_conti}
Let $h:\Cbar \to \Cbar$ be the semiconjugacy associated with a degeneration pair $(f \to g)$ that is given in \thmref{thm_semiconj}. Then $h$ tends to identity as $f$ tends to $g$.
\end{thm}

Here are two immediate corollaries: 
\begin{cor}\label{cor_conti_tiles}
The closures of $T_f(\theta,m,\ast)$ and $\Pi_f(\theta,\ast)$ in $\Tess(f)$ uniformly converge to those of $T_g(\theta,m,\ast)$ and $\Pi_g(\theta,\ast)$ in $\Tess(g)$ in the Hausdorff topology. 
\end{cor}

\begin{cor}\label{cor_diam_inv_arc_sys}
As $f \to g$, the diameters of connected components of $I_f$ uniformly tends to $0$.
\end{cor}

Let us start with some terminologies for the proof. Two degeneration pair $(f_1 \to g_1)$ and $(f_2 \to g_2)$ are \textit{equivalent} if $g_1=g_2$ and both $f_1$ and $f_2$ are in the same hyperbolic component. For a degeneration pair $(f \to g)$ by $f \approx g$ we mean $f$ is sufficiently close to $g$. In other words, the multiplier $r\omega$ of $O_f$ is sufficiently close to $\omega$, i.e., $r \approx 1$. 

Formally we consider a family of equivalent degeneration pairs $\skakko{(f \to g)}$ parameterized by $0<r<1$ and its behavior when $r$ tends to 1. To show the theorem, it suffices to show the following:
\begin{enumerate}[(i)]
\item For any compact set $K$ in $\Cbar- K_g$, we have $K \subset \Cbar-K_f$ for all $f \approx g$ and $h \to \id$ on $K$.
\item For any compact set $K$ in $\Kgc$, we have $K \subset \Kfc$ for all $f \approx g$ and $h \to \id$ on $K$.
\item $h$ is equicontinuous as $f \to g$ on the sphere. 
\end{enumerate}

In fact, any sequence $h_k$ associated with $f_k \to g$ has a subsequential limit $h_\infty$ which is identity on $\Cbar - J_g$ and continuous on $\Cbar$. Since $\Cbar-J_g$ is open and dense, the map $h_\infty$ must be identity on the whole sphere.

\subsection{Proof of (i)}
Let $\Boe_f:\Cbar-\D \to \Cbar-\Kfc$ be the extended B\"ottcher coordinate of $K_f$, i.e., $\Boe_f:\Cbar-\Dbar \to \Cbar-K_f$ is a conformal map with $\Boe_f(w^2)=f(\Boe_f(w))$; $\Boe_f(w)/w \to 1$ as $w \to \infty$; and $\Boe_f(e^{2 \pi i \theta}):=\gam_f(\theta) \in J_f$. Now (i) follows immediately from this stronger claim:

\begin{thm}[B\"ottcher convergence]\label{thm_Boettcher_converges}
As $f \to g$, we have a uniform convergence $\Boe_f \to \Boe_g$ on $\Cbar-\D$.
\end{thm}
Note that the uniform convergence on compact sets in $\Cbar-\Dbar$ is not difficult. Our proof is a mild generalization of the proof of Theorem 2.11 in \cite{Po}.

\begin{pf}
By \corref{cor_haus_conv_julia} one can easily check that $\Cbar-K_f$ converges to $\Cbar-K_g$ in the sense of Carath\'eodory kernel convergence with respect to $\infty$. Thus pointwise convergence $\Boe_f \to \Boe_g$ on each $z \in \Cbar-\Dbar$ is given by \cite[Theorem 1.8]{Po} and $\Boe_f'(\infty)=\Boe_g'(\infty)=1$. To show the theorem, it is enough to show that $K_f$ is uniformly locally connected as $f \to g$ by \cite[Corollary 2.4]{Po}. That is, for any $\e>0$, there exists $\delta>0$ such that for any $f \approx g$ and any $a,b \in K_f$ with $|a-b|<\delta$, there exists a continuum $E \subset K_f$ such that $a , b \in E$ and $\diam E < \e$. 

Suppose we have a sequence of equivalent degenerating pairs $(f_n \to g)$ such that: $f_n \to g$ uniformly; for $f_n$ there exist $a_n$ and $a_n'$ in $J_{f_n}$ with $|a_n -a_n'| \to 0$ which can not be contained in the same continuum in $K_{f_n}$ of diameter less than $\e_0>0$. We may set $a_n=\gam_{f_n}(\theta_n)$ and $a_n'=\gam_{f_n}(\theta_n')$ for some $\theta_n,~\theta_n' \in \T$ since $\gam_{f_n}$ maps $\T$ onto $J_{f_n}$. By passing to a subsequence, we may also assume that $\theta_n \to \theta$ and $\theta_n' \to \theta'$. Since $\gam_{f_n} \to \gam_g$ uniformly by \corref{cor_gam_f_to_gam_g}, the assumption $|a_n - a_n'| \to 0$ implies that we have $\gam_g(\theta)=\gam_g(\theta')=:w \in J_g$. 

\parag{Case 1:}
$\theta = \theta'$. We may assume that $\theta_n \le \theta_n'$ and both tend to $\theta$. Set $E_n:=\skakko{\gam_{f_n}(t):t \in [\theta_n, \theta_n']}$, which is a continuum containing $a_n$ and $a_n'$. Then for any $t \in [\theta_n, \theta_n']$, we have $|\gam_{f_n}(t)-w| \le |\gam_{f_n}(t)-\gam_g(t)|+|\gam_{g}(t)-\gam_g(\theta)| \to 0$ since $\gam_{f_n} \to \gam_g$ uniformly and $\gam_g$ is continuous. This implies $\diam E_n \to 0$ and is a contradiction. 

\parag{Case 2-1:}
$\theta \neq \theta'$ and $w \notin I_g$. First we show that $\gam_{f_n}(\theta)=\gam_{f_n}(\theta')$. Let $h_n:J_{f_n} \to J_g$ be the semiconjugacy given by \thmref{thm_semiconj_on_julia}. Since $h_n \cc \gam_{f_n}=\gam_g$, we have 
$$
w \ee h_n \cc \gam_{f_n}(\theta) = h_n \cc \gam_{f_n} (\theta') ~\notin~I_g.
$$
By property 1 of \thmref{thm_semiconj_on_julia}, this implies $\gam_{f_n}(\theta)=\gam_{f_n}(\theta')$. Now set
$$
E_n:=\skakko{\gam_{f_n}(t): |t-\theta|\le|\theta_n-\theta|
~\text{or}~|t-\theta'|\le|\theta_n'-\theta'|}, 
$$
which is a continuum containing $a_n$ and $a_n'$. Again one can easily check that $|\gam_{f_n}(t)-w| \to 0$ uniformly for any $\gam_{f_n}(t) \in E_n$ and thus $\diam E_n \to 0$.

\parag{Case 2-2:}
$\theta \neq \theta'$ and $w \in I_g$. There exists an $m \ge 0$ such that $g^m(w)=\beta_0$. Since $h_n \cc \gam_{f_n}=\gam_g$ we have $\gam_{f_n}(\theta), \gam_{f_n}(\theta') \in h_n^{-1}(w) \subset J_{f_n} \cap I_{f_n}$. If $q=1$, then $h_n$ is homeomorphism by \thmref{thm_semiconj_on_julia}. Thus $\gam_{f_n}(\theta)=\gam_{f_n}(\theta')$ and a contradiction follows from the same argument as above.

Suppose $q>1$. Then Case (a) ($q=q'$ and $l=l'$) by \propref{prop_case_a_case_b}. In particular, we have $w_n \in \al_{f_n}$ such that $w_n \to w$ and $f_n^m(w_n)$ is an attracting periodic point $\al_{0,n} \in O_{f_n}$ which tends to $\beta_0$. Let $\lam_n=r_n e^{2 \pi i p/q}$ be the multiplier of $O_{f_n}$ with $r_n \nearrow 1$. On a fixed small neighborhood of $w$, we have 
\begin{align*}
f^{-m} \cc f^{lq} \cc f^{m}(z) 
& \ee r_n^q z \, (1 + z^{q}+O(z^{2q})) \\
\longrightarrow~ 
g^{-m} \cc g^{lq} \cc g^{m}(z)
& \ee z  \, (1 + z^{q}+O(z^{2q}))
\end{align*}
by looking through suitable local coordinates as in Appendix A.2. (For simplicity, we abbreviate conjugations by the local coordinates.)

By \lemref{lem_small_inv_path}, we can find a small continuum $E_n' \subset K_{f_n}$ which joins $w_n$ and preperiodic points $\gam_{f_n}(\theta),~\gam_{f_n}(\theta')$. Set $E_n$ as in Case 2-1. Now $E_n' \cup E_n$ is a continuum containing $a_n$ and $a_n'$. Since $\diam (E_n' \cup E_n) \to 0$, we have a contradiction again.
\QED
\end{pf}

\subsection{Proof of (ii)}
Let us start with the following theorem:
\begin{thm}[Linearization convergence]\label{thm_linearization_converges}
Let $K$ be any compact set in $\Kgc$. Then $K \subset \Kfc$ for $f \approx g$ and $\Phi_f \to \Phi_g$ uniformly on $K$.
\end{thm}

\begin{pf}
One can easily check that $K \subset \Kfc$ if $f \approx g$ by \corref{cor_haus_conv_julia}. Let $\beta_0 \in O_g \cap \partial B_0'$. We may assume that $K'=g^N(K)$ is sufficiently close to $\beta_0$ and contained in $B_0'$ by taking a suitable $N \gg 0$. Then $K'$ is attracted to $\beta_0$ along the attracting direction associated with $B_0'$ by iteration of $g^{l'q'}$. For simplicity, set $\lbar:=lq=l'q'$.

Recall that $\Phi_f$ and $\Phi_g$ semiconjugate $f^\lbar$ and $g^\lbar$ to $F^q$ and $G^q$ in the fundamental model respectively. We will construct other semiconjugacies $\tilde{\Phi}_f$ and $\tilde{\Phi}_g$ with the same property as $\Phi_f$ and $\Phi_g$,  plus $\tilde{\Phi}_f \to \tilde{\Phi}_g$ on compact subsets of a small attracting petal in $B_0'$. Then we will show that they coincide.

By Appendix A.2, there exist local coordinates $\zeta=\psi_f(z)$ and $\zeta=\psi_g(z)$ with $\psi_f \to \psi_g$ near $\beta_0$ such that we can view $f^\lbar \to g^\lbar$ as 
\begin{align*}
f^\lbar(\zeta) 
& \ee \Lam \zeta \,(1 + \zeta^{q'}+O(\zeta^{2q'})) \\
\longrightarrow~ 
g^\lbar(\zeta) 
& \ee \zeta  \,(1+ \zeta^{q'}+O(\zeta^{2q'}))
\end{align*}
where $\Lam \to 1$. (To simplify notation, we abbreviate conjugations by these local coordinates. For example, by $f^\lbar(\zeta)$ we mean $ \psi_f \cc f^\lbar \cc \psi_f^{-1}(\zeta)$.) Now there are two cases for $\Lam$:
\begin{itemize}
\item In Case (a) ($q=q'$ and $l=l'$), the fixed point $\zeta=0$ is attracting and $\Lam=\lam^q=r^q=R<1$. 
\item In Case (b) ($q=1 < q' = l/l'$), the fixed point $\zeta=0$ is repelling and $|\Lam|>1$. 
\end{itemize}
Next by taking branched coordinate changes $w=\Psi_f(\zeta)=-\Lam^{q'}/(q'\zeta^{q'})$ and $w=\Psi_g(\zeta)=-1/(q'\zeta^{q'})$ respectively, we can view $f^\lbar \to g^\lbar$ as
\begin{align*}
f^\lbar(w) 
& \ee \Lam^{-q'} w  + 1 +O(1/w) \\
\longrightarrow~ 
g^\lbar(w) 
& \ee w  + 1 + O(1/w).
\end{align*}

\parag{Case (a).} 
Set $\tau=\Lam^{-q'}=R^{-q}>1$. By simultaneous linearization in Appendix A.3, we have convergent coordinate changes $W=u_f(w) \to u_g(w)$ on compact sets of $P_\rho:=\skakko{\rp w > \rho \gg 0}$ such that $f^\lbar \to g^\lbar$ is viewed as
\begin{align*}
\tilde{F}(W):=f^\lbar(W) 
& \ee \tau W  + 1 \\
\longrightarrow~ 
\tilde{G}(W):=g^\lbar(W) 
& \ee W  + 1.
\end{align*}
Let us adjust $\tilde{F} \to \tilde{G}$ to $F^q \to G^q$ in the fundamental model. Recall that the map $F(W)=RW+1$ has the attracting fixed point at $a=1/(1-R)$. On the other hand, the map $\tilde{F}$ has the repelling fixed point $\tilde{a}=1/(1-R^{-q})$ instead. Set $T_f(W):=aW/(W-\tilde{a})$. Then $T_f(W) = qW(1+O(W/\tilde{a})) \to T_g(W)=qW$ on any compact sets on the $W$-plane as $R \to 1$. By taking conjugations with $T_f$ and $T_g$, we can view $\tilde{F} \to \tilde{G}$ as $F^q \to G^q$ on any compact sets of the domain of $\tilde{G}$.

\parag{Case (b).} 
By Rouch\'e's theorem, there exists a fixed point $b$ of $f^\lbar(w) = \Lam^{-q'} w  + 1 +O(1/w)$ of the form $b=1/(1-\Lam^{-q'})+O(1)$. Indeed, this $b$ belongs to the image of the attracting cycle $O_f$ hence its multiplier is $r<1$. Set $S_f(w):=bw/(b-w)$. Then $S_f(w)=w(1+O(w/b)) \to S_g(w)=w$ on any compact sets of the $w$-plane as $r \to 1$. By taking conjugations by $S_f$ and $S_g$, we can view $f^\lbar(w) \to g^\lbar(w)$ as
\begin{align*}
f^\lbar(w) 
& \ee \tau w  + 1 +O(1/w) \\
\longrightarrow~ 
g^\lbar(w) 
& \ee w  + 1 + O(1/w)
\end{align*}
where $\tau=1/r>1$. By simultaneous linearization, we have convergent coordinate changes $W=u_f(w) \to u_g(w)$ on compact sets of $P_\rho$ such that $f^\lbar \to g^\lbar$ is again viewed as
\begin{align*}
\tilde{F}(W):=f^\lbar(W) 
& \ee \tau W  + 1 \\
\longrightarrow~ 
\tilde{G}(W):=g^\lbar(W) 
& \ee W  + 1.
\end{align*}
Since $q=1$, we adjust $\tilde{F} \to \tilde{G}$ to $F \to G$ in the fundamental model. Set $\tilde{b}:=1/(1-\tau)$ and $T_f(W):=\tilde{b}W/(\tilde{b}-W)$. Then $T_f(W) = W(1+O(W/\tilde{b})) \to T_g(W)=W$ on any compact sets on the $W$-plane as $r \to 1$. By taking conjugations by $T_f$ and $T_g$, we can view $\tilde{F} \to \tilde{G}$ as $F \to G$ on any compact sets of the domain of $\tilde{G}$.

\parag{Adjusting critical orbits.} 
Now we denote these final local coordinates conjugating $f^\lbar \to g^\lbar$ to $F^q \to G^q$ by $\hat{\Phi}_f \to \hat{\Phi}_g$, where the convergence holds on compact subsets of a small attracting petal $P'$ in $B_0'$ corresponding to $P_\rho$ in the $w$-plane.

We need to compare the images of the critical orbits by $\hat{\Phi}_f \to \hat{\Phi}_g$ on the $W$-plane and those by $\Phi_f$ and $\Phi_g$, and adjust their positions. 
We may assume that $g^{n\lbar}(0) \in P'$ for fixed $n \gg 0$. Then $f^{n\lbar}(0) \in P'$ for all $f \approx g$. Set $s:=\hat{\Phi}_f(f^{n\lbar}(0))$ and $s':=\hat{\Phi}_g(g^{n\lbar}(0))$. Then $s \to s'$ as $f \to g$. On the other hand, we have 
$$
\Phi_f(f^{n\lbar}(0))= F^{nq}(\Phi_f(0)) = F^{nq}(0) = R^{nq-1}+\cdots +1=:R_n
$$
and $\Phi_g(g^{n\lbar}(0))=nq$. Set $U_f(W):=k(W-a)+a$ and $U_g(W):=W+nq-s'$ where $k=(R_n-a)/(s-a)$. Then one can check that $U_f \to U_g$ on any compact sets in the $W$-plane as $f \to g$ and $U_f$ and $U_g$ commute with $F$ and $G$ respectively. By defining $\tilde{\Phi}_f$ and $\tilde{\Phi}_g$ by $U_f \cc \hat{\Phi}_f$ and $U_g \cc \hat{\Phi}_g$ respectively, we have $\tilde{\Phi}_f \to \tilde{\Phi}_g$ on compact sets of $P'$ with $\tilde{\Phi}_f(f^{n\lbar}(0))=R_n$ and $\tilde{\Phi}_g(g^{n\lbar}(0))=nq$.

Finally we need to check that $\tilde{\Phi}_f=\Phi_f$ and $\tilde{\Phi}_g=\Phi_g$. The latter equality is clear by uniqueness of the Fatou coordinate (\cite[\S 8]{MiBook}). For the former, recall that $W=\Phi_f(z)$ is given by
$$
z ~\mapsto~ \phi_f(z)=w ~\mapsto~ w^q =W~\mapsto~ a(1-W)=:\Phi_f(z)
$$
and $\phi_f$ is uniquely determined under the condition of $\phi_f(0)=1$ (\cite[\S 10]{MiBook}). Let us consider the local coordinate $\tilde{\phi}_f$ on a compact set of $P'$ given by 
$$
z ~\mapsto~ \tilde{\Phi}_f(z)=W ~\mapsto~ 
\kakko{1-\frac{W}{a}}^{1/q}=:\tilde{\phi}_f(z)=w,
$$
where we take a suitable branch of $q$th root such that $\tilde{\phi}_f(f^{n\lbar}(0))=\lam^{nq}$ on the $w$-plane. Then $\tilde{\phi}_f(f(z))=\lam \tilde{\phi}_f(z)$. Since $\phi_f(0)=1$ is equivalent to $\phi_f(f^{n\lbar}(0))=\lam^{nq}$, the map $\tilde{\phi}_f$ coincide with $\phi_f$. This implies the former equality. 

Now we may assume that $K'=g^{N}(K) \subset D \Subset P'$ for some open set $D$. If $f \approx g$, then $f^N(K) \subset D$ and we have the uniform convergence $\Phi_g \to \Phi_f$ on $D$. We finally obtain the uniform convergence on $K$ by $\Phi_f(z)=F^{-N}(\Phi_f(f^N(z))) \to G^{-N}(\Phi_g(g^N(z))=\Phi_g(z)$ for $z \in K$.
\QED
\end{pf}

\paragraph{Proof of (ii).} 
We first work with the fundamental model. Suppose $\e \searrow 0$, and set $R=1-\e$. Then $F(W)=RW +1$ fixes $a_\e=1/(1-R)=\e^{-1}$. For a fixed $\gam$ with $1/2<\gam<1$, we define a compact set $Q_\e \subset \C$ by:
$$
Q_\e \dee \skakko{W=a_\e+\rho e^{(\pi-t)i} \in \C~:~
|t| \le \e^\gam,~|\rho-a_\e| \le a_\e \sin \e^\gam
}.
$$
Let $D$ be any bounded set in $\C$. For all $\e \ll 1$, the compact set $Q_\e$ contains $D$. Let $H:\C-[a_\e, \infty) \to \C$ be the conjugacy between $F$ and $G(W)=W+1$ as in \secref{sec_04}. Then one can easily check that $|\rp W-\rp H(W)|=O(\e^{2\gam-1})$ and $|\ip W-\ip H(W)|=O(\e^{2\gam-1})$ on $Q_\e$. Thus $H \to \id$ uniformly on $D$. 

Let $K$ be any compact set in $\Kgc$, and let $D$ be the 1/10-neighborhood of $\Phi_g(K)$. For all $f \approx g$, we have $K \subset \Kfc$ and $\Phi_f(K) \subset D$ by \thmref{thm_linearization_converges}. By the argument above on the fundamental model, the restriction $h|_K$ is a branch of $\Phi_g^{-1} \cc H \cc \Phi_f$ that converges to the identity. (The branch is determined by the tile-to-tile correspondence given by $h$.) \QED

\subsection{Proof of (iii)}
To show (iii) we need two propositions on properties of panels as $f \to g$. The first one is a refinement of \propref{prop_diam}, and the second one is on the convergence of panels with a fixed angle: 

\begin{prop}[Uniformly small panels]
\label{prop_uniformly_small_panels}
For any $\e>0$, there exists $N=N(\e)$ such that for all $f \approx g$, $\ast=\pm$ and $\theta \in \Theta_g$ with $\depth(\theta) \ge N$, 
$$
\diam \Pi_f(\theta, \ast)<\e ~~~\textit{and}~~~
\diam \Pi_g(\theta, \ast)<\e.
$$
\end{prop}

\begin{prop}[Hausdorff convergence to a panel]
\label{prop_panel_conv}
For fixed angle $\theta \in \Theta_g$ and signature $\ast=+$ or $-$, we have $\overline{\Pi_{f}(\theta,\ast)} \to \overline{\Pi_{g}(\theta,\ast)}$ as $f \to g$ in the Hausdorff topology. 
\end{prop}

Let us show (iii) first by assuming them:

\parag{Proof of (iii).}
By (i) we have equicontinuity near $\infty$. Assume that there exist degeneration pairs $(f_k \to g)$ with semiconjugacies $h_k$ as in \thmref{thm_semiconj}, points $a_k,~a_k' \in \C$ with $|a_k-a_k'| \to 0$, and $b_k=h_k(a_k),~b_k'=h_k(a_k')$ with $|b_k-b_k'| \ge \e_0>0$. 

Suppose that $a_k,~a_k' \in \C-\Kfc$ thus $b_k,~b_k' \in \C-\Kgc$. Then there exists $w_k,~w_k' \in \C-\D$ such that $a_k=B_{f_k}(w_k),~a_k'=B_{f_k}(w_k')$ and $b_k=B_{g}(w_k),~b_k'=B_{g}(w_k')$. By \thmref{thm_Boettcher_converges}, we have $B_{f_k} \to B_{g}$. Thus $|a_k-a_k'| \to 0$ implies $|b_k-b_k'| \to 0$, a contradiction.

Now it is enough to show the case where $a_k,~a_k' \in K_{f_k}$ thus $b_k,~b_k' \in K_g$. By taking subsequences, we may assume that $a_k \to a$, $a_k' \to a$, $b_k \to b$, and $b_k' \to b'$ with $|b-b'| \ge \e_0/2>0$. Since $K_{f_k} \to K_g$ in the Hausdorff topology, $a, ~b$ and $b'$ are all in $K_g$. 

First let us consider the case where $a$ is bounded distance away from $J_g$. Then we have a compact neighborhood $E$ of $a$ such that $h_k|_E \to \id|_E$ and $a_k, a_k' \in E$ for all $k \gg 0$. This implies that $|b_k- b_k'| \to 0$, a contradiction.

Next we consider the case where $a \in J_g$. For $a_k \to a$ and $b_k \to b$, we will claim that $a=b$. Then by the same argument we have $a=b'$ and this is a contradiction. 

For $a_k \in K_{f_k}$, we take any $\theta_k \in \T$ such that: $a_k=\gam_{f_k}(\theta_k)$ if $a_k \in J_{f_k}$; otherwise $a_k$ is contained in the closure of $\Pi_{f_k}(\theta_k)$. (Then $b_k=\gam_g(\theta_k)$ or $b_k$ is in the closure of $\Pi_{g}(\theta_k)$.) By passing to a subsequence, we may assume that $\theta_k \to \theta$ for some $\theta \in \T$. 

If $\theta_k \notin \Theta_g$, we define its depth by $\infty$. Then there are two more cases according to $\limsup \depth(\theta_k) =\infty$ or not.

If $\limsup \depth(\theta_k) =\infty$, we take a subsequence again and assume that $\depth(\theta_k)$ is strictly increasing. Then by \propref{prop_uniformly_small_panels} we have $|a_k-\gam_{f_k}(\theta_k)| \to 0$. Since $\theta_k \to \theta$ and $\gam_{f_k} \to \gam_{g}$ uniformly (\corref{cor_gam_f_to_gam_g}), we have $|a_k -\gam_g(\theta)| \to 0$, thus $a = \gam_g(\theta)$. Similarly we conclude that $b=\gam_g(\theta)$ and this implies a contradiction. 

If $\limsup \depth(\theta_k) < \infty$, we take a subsequence again and assume that $\theta_k=\theta \in \Theta_g$ for all $k \gg 0$. By \propref{prop_panel_conv} $a_k \in \overline{\Pi_{f_k}(\theta)}$ are approximated by some $c_k \in \Pi_g(\theta)$ with $|a_k - c_k| \to 0$ thus $c_k \to a \in J_g$. Since $\overline{\Pi_g(\theta)} \cap J_g=\skakko{\gam_g(\theta)}$, we have $a=\gam_g(\theta)$. On the other hand, if $b_k \in \Pi_g(\theta)$ is bounded distance away from $J_g$, there exists a compact neighborhood $E' \subset \Kgc$ of $b$ where $h_k|_{E'} \to \id|_{E'}$ and it leads to a contradiction. Thus $b \in J_g$ and it must be $\gam_g(\theta)$. Now we obtain $a=b$. \QED

\paragraph{}
To complete the proof of \thmref{thm_conti}, let us finish the proofs of the propositions.

\parag{Proof of \propref{prop_uniformly_small_panels}.}
We modify the argument of \propref{prop_diam}. Suppose that there exist $f_k \to g$ which determine equivalent degeneration pairs $(f_k \to g)$ and $\theta_k$ with $n_k=\depth(\theta_k) \nearrow \infty$ such that $\diam \Pi_{f_k}(\theta_k, +) \ge \e_0>0$ for all $k$. Then we can take a branch $F_k$ of $f_k^{-n_k}$ such that $F_k$ maps $\Pi_{f_k}(\theta_0^+, +)^\cc$ onto $\Pi_{f_k}(\theta_k, +)^\cc$ univalently. 

Take a small ball $B \Subset T_g(\theta_0^+,0,+)$ and fix a point $z \in B$. By (ii), we may assume that $B \Subset T_{f_k}(\theta_0^+,0,+)$ for all $k \gg 0$. Since $F_k|_B$ avoid values near $\infty$, they form a normal family. By passing to a subsequence, we may also assume that there exists $\phi=\lim F_k|_B$ that is non-constant by assumption. Now we have a small open set $V \Subset \phi(B)$ with $V \subset F_k(B)$ for all $k \gg 0$, thus $f_k^{n_k}(V) \subset B \subset K_{f_k}^\cc$. This implies that $V \subset K_{f_k}^\cc$ for all $k \gg 0$ hence by \corref{cor_haus_conv_julia} we have $V \subset \Kgc$ too. Since $V$ is open, there exist a tile $T=T_g(\theta, m,+)$ and a small ball $B'$ such that $B' \Subset (T \cap V)^\cc$. By $B' \subset T$ and (ii) again, we have $B' \subset T_k:=T_{f_k}(\theta, m, +)$ for all $k \gg 0$. Moreover, since $B' \subset V$ we have $f_k^{n_k}(B') \subset B \subset T_{f_k}(\theta_0^+,0,+)$. Thus $f_k^{n_k}(T_k)$ must be $T_{f_k}(\theta_0^+,0,+)$ but $f_k^{n_k}(T_k)$ has level $m+n_k \to \infty$. It is a contradiction. 

Finally one can finish the proof by the same argument as \propref{prop_diam}.
\QED

\paragraph{Proof of \propref{prop_panel_conv}.}
It is enough to consider the case of $\theta=\theta_0^+$ and $\ast=+$. Recall that the attracting cycle $O_f$ has the multiplier $re^{2 \pi i p/q}$. We introduce a parameter $\e \in [0,1)$ of $f \to g$ such that $r^q=R=1-\e$. Set $\Pi_\e:=\Pi_f(\theta_0^+, +)$ and $\Pi_0:=\Pi_g(\theta_0^+, +)$. Then the semiconjugacy $h=h_\e$ sends $\Pi_\e$ to $\Pi_0$. To conclude the statement it is enough to show the following: For any $\delta>0$, we have $\Pi_0 \subset N_\delta(\Pi_\e)$ and $\Pi_\e \subset N_\delta(\Pi_0)$ for all $\e \ll 1$, where $N_\delta(\cdot)$ denotes the $\delta$-neighborhood.

It is easy to check $\Pi_0 \subset N_\delta(\Pi_\e)$: We can take a compact set $K$ such that $K \subset \Pi_0^\cc \Subset N_\delta(K)$. Since $h_\e \to \id$ on $K$, we have $K \subset \Pi_\e^\cc$ for all $\e \ll 1$. Thus we have $\Pi_0 \subset N_\delta(K) \subset N_\delta(\Pi_\e)$.

The proof of $\Pi_\e \subset N_\delta(\Pi_0)$ is more technical. Here let us assume that $q=q'$, Case (a). Case (b) ($q=1<q'$) is merely analogous and left to the reader.

\parag{Local coordinates.}
Set $B:=B(\beta_0, \delta)$. For fixed $\delta$ that is small enough, there exists a convergent family of local coordinates $\zeta=\psi_\e(z) \to \psi_0(z)$ on $B$ with the following properties for all $0 \le \e \ll 1$: 
\begin{itemize}
\item There exists $\delta' >0$ 
independent of $0 \le \e \ll 1$ such that $\Delta:=B(0,\delta') \Subset \psi_\e(B)$.
\item Set $\fe:=f^{lq}$, $f_0:=g^{lq}$, and $R_\e:=1-\e$. Then $\fe(\zeta)=R_\e \zeta \, (1+\zeta^q+O(\zeta^{2q}))$ on $\Delta$. (See Appendix A.2.)
\item $\psi_\e$ maps $\Pi_\e \cap \psi_\e^{-1}(\Delta)$ into $\Delta':=\skakko{\zeta \in \Delta: -\pi/2q < \arg \zeta < 3 \pi/2q }$. (This is just a technical assumption.)
\item Set $E_\e:=\skakko{\zeta \in \Delta': |\arg \zeta^q| \le \pi/3,~|\zeta^q| \ge \e/2}$. Then $f_0^{-1}(\overline{E_0}) \subset E_0 \cup \skakko{0}$ and $\fe^{-1}(\overline{E_\e}) \subset E_\e$ for all $0< \e \ll 1$. (See the argument of \lemref{lem_small_inv_path}). 
\end{itemize}

Let us interpret the setting of \thmref{thm_linearization_converges} by using $\e \in [0,1)$. If $0 < \e < 1$, we denote $\Phi_f$, $\Psi_f$, $u_f$, $T_f$, and $U_f$ by $\Phi_\e$, $\Psi_\e$, $u_\e$, $T_\e$, and $U_\e$ respectively. If $\e=0$ they denote $\Phi_g$, $\Psi_g$, etc. In particular, we consider $\Psi_\e$ only on $\Delta'$. For later use, we define $W=\chi_\e(\zeta)$ for each $\zeta \in \psi_\e(K_{\fe}^\cc \cap B)$ by $\chi_\e:=\Phi_\e \cc \psi_\e^{-1}$. 

Now we can see $\fe$ on $\Delta'$ through $w=\Psi_\e(\zeta)$ as $\fe(w)=\taue w+1+O(1/w)$ where $\taue:=R_\e^{-q}$. On this $w$-plane, take $P=P_\rho=\skakko{\rp w \ge \rho \gg 0}$ such that for all $0 \le \e \ll1$, the set $\hat{P}:=\Psi_\e^{-1}(P)$ is contained in $\Delta'$ and that $u_\e$ is defined on $P$. Note that for all $0 \le \e \ll 1$ we have $\fe(P) \subset P$ and $u_0(w)=w(1+o(1))$ by \lemref{lem_u_0}. One can also check that $\chi_\e \cc \Psi_\e^{-1}(w)=U_\e \cc T_\e \cc u_\e(w)$ on $P$ and
$$
 U_\e \cc T_\e \cc u_\e(w) \ee 
 U_0 \cc T_0 \cc u_0(w) + o(1) \ee qw(1+o(1))
$$ 
on compact sets of $P$.

\parag{Rectangles.}
For fixed positive integers $M$ and $N$, we define the following compact sets in the $W$-plane:
\begin{align*}
C_0 &\dee \skakko{W \in \C ~:~ 
(N-1)q \le \rp W \le Nq, ~ 0 \le \ip W \le Nq} \\
Q_0&\dee \bigcup_{k=0}^M G^{-kq}(C_0)
~~~\text{and}~~~C_0' \dee G^{-Mq}(C_0)
\end{align*}
where $G(W)=W+1$. 

\begin{figure}[htbp]
\centering{
\includegraphics[width=0.8\textwidth]{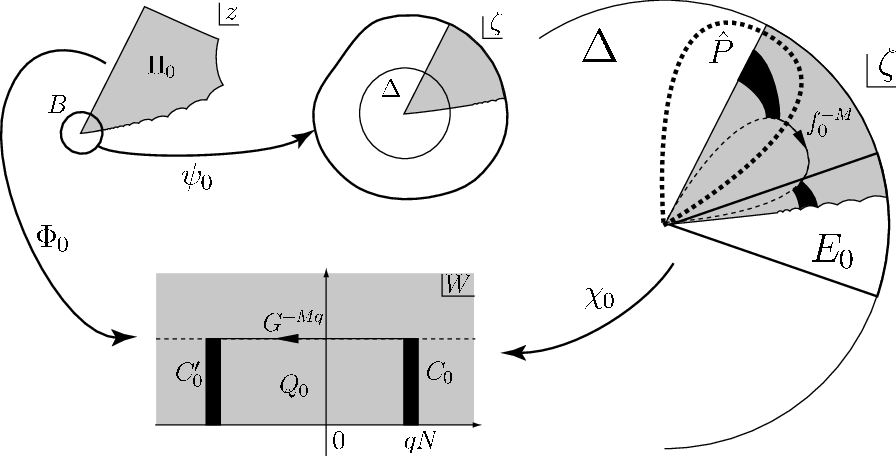}
}
\caption{Taking $M$ and $N$.}\label{fig_coordinates}
\end{figure}

By taking sufficiently large $N$ and $M$, we may assume the following:
\begin{enumerate}[(1)]
\item Set $\tilde{Q}_0:=\Pi_0 \cap \Phi_0^{-1}(Q_0)$ in the $z$-coordinate. Then $\Pi_0-\tilde{Q}_0 \Subset \psi_0^{-1}(\Delta)$.
\item In the $w$-coordinate, we have $\chi_0^{-1}(C_0) \subset \hat{P}$ and $\chi_0^{-1}(C_0') \subset E_0$.
\end{enumerate}
See \figref{fig_coordinates}. In fact, for any compact set $K$ with $\Pi_0-K \Subset \psi_0^{-1}(\Delta)$, the set $\Phi_0(K)$ is compact in $\overline{\Hyp}_W:=\skakko{\ip W \ge 0}$ and covered by $Q_0$ if we take sufficiently large $N$ and $M$. Thus we have (1). If $N \gg 0$, the set $C_0$ must be contained in $\chi_0(P)$. Since $\Pi_0 \cap \Phi_0^{-1}(C_0)$ is compact, it is uniformly attracted to the repelling direction by iteration of $(g|_{\Pi_0})^{-lq}$. Thus we have (2) by taking $M$ much larger.

\parag{Perturbation.} 
We fix such integers $N$ and $M$. Now we consider perturbation of fixed rectangles $C_0, C_0'$, and $Q_0$ with properties (1) and (2). By using the conjugacy $H=H_\e: \C-[a, \infty) \to \C$ between $F=F_\e$ and $G=F_0$, we define $C_\e$, $C_\e'$ and $Q_\e$ by their homeomorphic images by $H_\e^{-1}$. Since $\He \to \id$ as $\e \to 0$ on any compact sets (see the proof of (ii)), we have $C_\e \to C_0$, $C_\e' \to C_0'$, and $Q_\e \to Q_0$ in the Hausdorff topology. Moreover, we have the following properties for all $\e \ll 1$:
\begin{enumerate}[(1')]
\item Set $\tilde{Q}_\e:=\Pi_\e \cap \Phi_\e^{-1}(Q_\e)$ in the $z$-coordinate. Then $\tilde{Q}_\e \subset N_{\delta/2}(\tilde{Q_0})$.
\item In the $\zeta$-coordinate, we have $\chi_\e^{-1}(C_\e) \subset \hat{P}$ and $\chi_\e^{-1}(C_\e') \subset E_\e$.
\end{enumerate}
In fact, since $\tilde{Q}_0 =h_\e(\tilde{Q}_\e)$ and is compact, property (1') follows by $\Phi_\e \to \Phi_0$ as $\e \to 0$. Property (2') holds because $\chi_\e \to \chi_0$ on compact sets in $\hat{P}$ and $f^{lqM} \to g^{lqM}$. 

Now it is enough to show $\Pi_\e-\tilde{Q}_\e \Subset \psi_\e^{-1}(\Delta) \subset B$, which is equivalent to $\chi_\e^{-1}(\overline{\Hyp}_W-Q_\e) \Subset \Delta$ in the $\zeta$-coordinate. We consider the following three sets in $\overline{\Hyp}_W$:
\begin{align*}
X_0 &\dee 
\skakko{W \in \overline{\Hyp}_W: \rp W \le (N-M-1)q,~\ip W \le Nq } \\
Y_0 &\dee 
\skakko{W \in \overline{\Hyp}_W: \rp W \ge Nq,~\ip W \le Nq} \\
Z_0 &\dee 
\skakko{W \in \overline{\Hyp}_W: \ip W \ge Nq} 
\end{align*}
Let $X_\e, Y_\e$, and $Z_\e$ be their homeomorphic images by $H_\e^{-1}$. Then $X_\e \cup Y_\e \cup Z_\e=\overline{\Hyp}_W-Q_\e^\cc$.

Note that $X_\e \ee \bigcup_{k \ge 1} F_\e^{-kq}(C_\e')$ and $Y_\e \ee \bigcup_{k \ge 1} F_\e^{kq}(C_\e)$. Since $\fe^{-1}(E_\e) \subset E_\e$ and $\fe(\hat{P}) \subset \hat{P}$ in the $\zeta$-coordinate, (2') implies $\chi_\e^{-1}(X_\e) \subset E_\e$ and $\chi_\e^{-1}(Y_\e) \subset \hat{P}$ thus we have $\chi_\e^{-1}(X_\e \cup Y_\e) \subset \Delta$.

The proof is completed by showing $\chi_\e^{-1}(Z_\e) \subset \Delta$. It is enough to show that $\chi_\e^{-1}(\partial Z_\e) \subset \Delta$. Note that $\partial Z_\e$ consists of two half lines, one is the interval $I_\e:=[a_\e, \infty)$ where $a_\e$ is the attracting fixed point of $F_\e$, and the other is $I_\e':=H_\e^{-1}(\partial Z_0)$, the one along the top edge of $Q_\e$. 

First we show that $\chi_\e^{-1}(I_\e) \subset \Delta$. Recall that this is the image of a degenerating arc in the $\zeta$-coordinate. Let $E_0':=\skakko{\zeta \in \Delta': |\arg (-\zeta^q)| \le \pi/3}$. Then one can check that $\fe(E_0') \subset E_0'$ and $\fe^{-1}(E_0) \subset E_0$ for all $\e \ll 1$ as in the argument of \lemref{lem_small_inv_path}. 

The real part of $g^{lqk}(0)$ in the $w$-coordinate increases as $k \to \infty$ thus the critical orbit of $f_0=g^{lq}$ in $\Delta'$ is tangent to the attracting direction in the $\zeta$-coordinate. Thus we may assume that $g^{lqn}(0)$ in the proof of \thmref{thm_linearization_converges} is contained in $E_0'$. Hence $f^{lqn}(0)=\fe^n(0)$ in the $\zeta$-coordinate is contained in $E_0'$ for all $\e \ll 1$. Moreover, the property $\fe(E_0') \subset E_0'$ implies that the critical orbit of $\fe=f^{lq}$ in $\Delta'$ is eventually contained in $E_0'$. By construction of the degenerating arcs in \lemref{lem_landing} and by $\fe^{-1}(E_0) \subset E_0$, the arc $\chi_\e^{-1}(I_\e)$ must be contained in $E_0 \subset \Delta$.

Next we show that $\chi_\e^{-1}(I_\e') \subset \Delta$. Let $s_\e$ and $\ell_\e$ be the top edges of quadrilaterals $C_\e$ and $Q_\e$ intersecting $I_\e'$. Then $\ell_\e=\bigcup_{k \ge 0}^M F_\e^{-kq}(s_\e)$. Now it is enough to show that $\chi_\e^{-1}(\ell_\e)$ is contained in $\Delta$ since $\chi_\e^{-1}(X_\e \cup Y_\e) \subset \Delta$. 

Take any point $w_0$ in $\Psi_\e \cc \chi_\e^{-1}(s_\e)=(U_\e \cc T_\e \cc u_\e)^{-1}(s_\e)$ in the $w$-plane. We may assume that $N$ is sufficiently large and $w_0 \in B(N+Ni, N/4)$ for all $\e \ll 1$, since $U_\e \cc T_\e \cc u_\e(w)=qw(1+o(1))$ on compact sets of $P$. Moreover, we may assume that $\Psi_\e(\partial \Delta) \subset B(0, N/4)$.

Recall that $\fe(w)=\taue w+1+O(1/w)$ and thus $\fe^{-1}(w)=\taue^{-1} (w-1)+O(1/w)$. Take any $w$ with $N/4 \le |w| \le 4N$. Then we have $|\fe^{-1}(w)-(w-1)|=O(\e N)+O(1/N)$. Thus for any fixed $\kappa \ll 1$, by taking $N \gg 0$ we have $|\fe^{-1}(w)-(w-1)| \le \kappa$ for all $\e \ll 1$. This implies $|\arg(\fe^{-1}(w)-w)| \lee \arcsin \kappa$.

\begin{figure}[htbp]
\centering{\vspace{0cm}
\includegraphics[height=4cm]{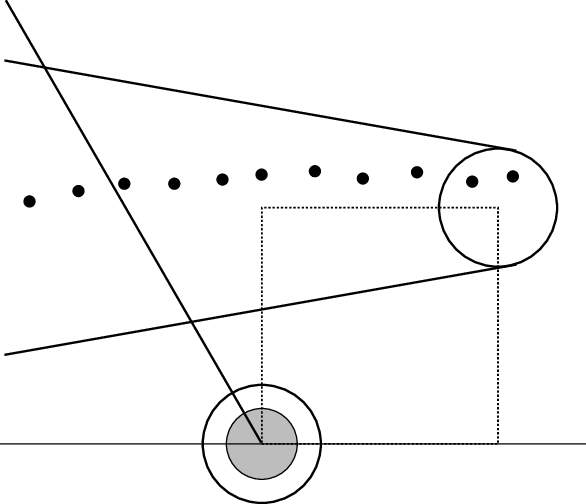}}
\caption{An orbit in the $w$-plane. The dotted square has height $N$.}\label{fig_w-plane}
\end{figure}
By (2'), the orbit $w_k=\fe^{-k}(w_0)$ of $w_0$ lands on $\Psi_\e(E_\e)$ by at most $M$ iteration of $\fe^{-1}$ (thus $\arg w_M >2 \pi/3 $). For small enough $\kappa$, the point $w_k$ satisfies $N/4 \le |w_k| \le 4N$ and $|\arg(w_k-w_0)| \lee \arcsin \kappa$ for all $k=0, \ldots, M$. (See \figref{fig_w-plane}.) This implies that $\Psi_\e \cc \chi_\e^{-1}(\ell_\e)$ never crosses over $\Psi_\e(\partial \Delta)$ thus we have $\chi_\e^{-1}(\ell_\e) \subset \Delta$. \QED

\setcounter{section}{0}
\renewcommand{\thesection}{\Alph{section}}

\section{Appendix}\label{sec_06}
In this section we give some investigation on perturbation of parabolic cycle corresponding to the degeneration pair $(f \to g)$.

\subsection{Pinching semiconjugacy on the Julia sets}
Let $(f \to g)$ be a general degeneration pair. Recall that the attracting cycle $O_f=\skakko{\llist{\al}=a_0}$ has multiplier $\lam=r\omega=r\exp(2 \pi i p/q)$ with $0<r<1$; and that the parabolic cycle $O_g=\skakko{\beta_1, \ldots, \beta_{l'}=\beta_0}$ has multiplier $\omega'=\exp(2 \pi p'/q')$. 

By applying Theorem 1.1 of \cite{Ka1} to $(f \to g)$, we have:

\begin{thm}\label{thm_semiconj_on_julia}
If $f \approx g$, there exists a unique semiconjugacy $h_J:J_f \to J_g$ with the following properties: 
\begin{enumerate}
\item If $\card \, h_J^{-1}(w) \ge 2$ for some $w \in J_g$ then $w \in I_g$ and $\card \, h_J^{-1}(w)=q$ (thus $q =q' \ge 2)$. 
\item $h_J$ is a homeomorphism iff $(f \to g)$ is of type $q=1$. 
\item $\displaystyle \sup_{z \in J_f} |z-h_J(z)| \to 0$ as $f \to g$.
\end{enumerate}
\end{thm}
(See also \propref{prop_case_a_case_b}.) The proof of Theorem 1.1 of \cite{Ka1} is based on a pull-back argument and it does not use quasiconformal maps. Here is a useful corollary which easily follows from property 3:
\begin{cor}\label{cor_haus_conv_julia}
As $f \to g$, the Julia set $J_f$ converges to $J_g$ in the Hausdorff topology.
\end{cor}

Since $h_J \cc \gam_f$ and $\gam_g$ determines the same ray equivalence, we have $h_J \cc \gam_f=\gam_g$. For $\theta \in \T$, put $\gam_f(\theta)$ into $z$ in property 3 of the theorem above. Then we have:

\begin{cor}\label{cor_gam_f_to_gam_g}
As $f \to g$, the map $\gam_f:\T \to J_f$ converges uniformly to $\gam_g:\T \to J_g$.
\end{cor}

\subsection{Normalized form of local perturbation}
For a degeneration pair $(f \to g)$, the parabolic cycle $O_g$ is approximated by an attracting or repelling cycle $O_f'$ with the same period $l'$ and multiplier $\lam' \approx \omega'=e^{2 \pi i p'/q'}$. (See \secref{sec_02}. 
Note that in the quadratic family, 
a parabolic cycle with multiplier $\omega'$ has $q'$ petals
if and only if $\omega'$ is a primitive $q'$th root of unity.)
Let $\al_0'\in O_f'$ with $\al_0' \to \beta_0$. Then by looking through the local coordinates $\psi_f(z)=z-\al_0'$ and $\psi_g(z)=z-\beta_0$ near $\beta_0$ one 
observe the convergence $f^{l'} \to g^{l'}$ as follows:
\begin{align*}
\psi_f \cc f^{l'} \cc \psi_f^{-1}(z) 
& \ee \lam' z +O(z^2)\\
\longrightarrow~ 
\psi_g \cc g^{l'} \cc \psi_g^{-1}(z) 
& \ee \omega' z +O(z^2).
\end{align*}
Here we claim that by replacing $\psi_f \to \psi_g$ with better local coordinates, we have a normalized form of convergence:
\begin{align*}
\psi_f \cc f^{l'} \cc \psi_f^{-1}(z) 
& \ee \lam' z  + z^{q'+1}+O(z^{2q'+1}) \\
\longrightarrow~ 
\psi_g \cc g^{l'} \cc \psi_g^{-1}(z) 
& \ee \omega' z  + z^{q'+1}+O(z^{2q'+1}).
\end{align*}
\newcommand{\lame}{\lambda_\epsilon}

More generally, we have:
\begin{prop}\label{prop_normal_form}
For $\e \in [0,1]$, let $\skakko{\fe}$ be a family of holomorphic maps on a neighborhood of $0$ such that as $\e \to 0$,
$$
\fe(z) = \lame z + O(z^2) 
~\longrightarrow ~
f_0(z) = \lam_0 z + O(z^2)
$$
where $\lam_0$ is a primitive $q$th root of unity
and $f_0(0)=0$ has $q$ petals. Then we have a family of holomorphic maps $\skakko{\phie}$ 
for $0 \le \e \ll 1$ such that
$$
\phie \cc \fe \cc \phie^{-1}(z) \ee \lame z + z^{q+1}+O(z^{2q+1})
$$
and $\phie \to \phi_0$ near $z=0$.
\end{prop}

\begin{pf}
\footnote{In the published version of this paper 
({\it Erg. Th. Dyn. Sys.} {\bf 29} (2009), 579--612),
``$\lame^{n}-\lame$" is incorrectly written as 
``$\lame^{n+1}-\lame$" in this proof.  
}
First suppose that $\fe(z)=\lame z + A_\e z^n+O(z^{n+1})$ where $2 \le n \le q$. Let us consider a coordinate change by $z \mapsto z - B_\e z^n$ with $B_\e=A_\e/(\lame^{n}-\lame)$. Note that 
$\lame^{n}-\lame$ is bounded distance away from $0$ when $\e \ll 1$, because $\lame$ converges to a primitive $q$th root of unity. In particular, the coordinate change $z \mapsto z - B_\e z^n$ also converges to $z \mapsto z - B_0 z^n$ near $0$. By applying these coordinate changes, we can view the family $\skakko{\fe}$ as
$$
\fe(z) \ee \lame z + O(z^{n+1}).
$$
By repeating this process until $n=q$, we have the family $\skakko{\fe}$ of the form 
$$
\fe(z) \ee \lame z + C_\e z^{q+1}+A_\e' z^n+O(z^{n+1})
$$
where $q+2\le n \le 2q$,
and we have $C_\e \neq 0$ for $0 \le \e \ll 1$
since $C_0 \neq 0$ 
(otherwise $f_0^q(z)=z+O(z^{q+2})$
and it contradicts the assumption that 
$f_0(0)=0$ has $q$ petals).
Next for each $0 \le \e \ll 1$
 take a linear coordinate change by $z \mapsto C_\e^{1/q} z$ to normalize $C_\e$ to be $1$.
By taking another coordinate change of the form $z \mapsto \zeta=z - B_\e' z^n$ with $B_\e'=A_\e'/(\lame^{n}-\lame)$ again, we have 
$$
\fe(z) \ee \lame z + z^{q+1}+O(z^{n+1}).
$$
We can repeat this process until $n=2q$ and we have the desired form of convergence. $\blacksquare$
\end{pf}

For this new family $\skakko{\fe(z)= \lame z + z^{q+1}+O(z^{2q+1})}$ and $n \ge 0$, one can easily check that 
$$
\fe^n(z) \ee \lame^n z + C_{\e,n} z^{q+1}+O(z^{2q+1})
$$
where $C_{\e,n}$ is given by the recursive formula $C_{\e, n+1}=\lame^{q+1} C_{\e, n} +\lame^n$. Let $n=q$ and set $\varLambda_\e:=\lame^q ~(\to 1$ as $\e \to 0$). By taking linear coordinate changes with $z \mapsto (C_{\e, q}/\varLambda_\e)^{1/q}z$, we have the convergence of the form
\begin{align*}
\fe^q(z) 
& \ee \varLambda_\e z \, (1 + z^{q}+O(z^{2q})) \\
\longrightarrow~ 
f_0^q(z) 
& \ee z \, (1 + z^{q}+O(z^{2q})).
\end{align*}
By further coordinate changes with $w=\Psi_\e(z)=-\varLambda_\e^q/(qz^q)$, we have 
\begin{align*}
\Psi_\e \cc \fe^q \cc \Psi_\e^{-1}(w) 
& \ee \varLambda_\e^{-q} w + 1 + O(1/w) \\
\longrightarrow~ 
\Psi_0 \cc f_0^q \cc \Psi_0^{-1}(w)
& \ee  w + 1 +O(1/w)
\end{align*}
on a neighborhood of infinity. Note that we have a similar representation for $f^{l'q'} \to g^{l'q'}$.

\subsection{Simultaneous linearization}
Recently T.Ueda \cite{Ue} showed the simultaneous linearization theorem that explains hyperbolic-to-parabolic degenerations of linearizing coordinates. Here we state a simple version of the theorem which is enough for our investigation. For $R \ge 0$, let $E_R$ denote the region $\skakko{z \in \C: \rp z  \ge R}$. 

\begin{thm}[Ueda]\label{thm_ueda}
For $\e \in [0,1]$, let $\skakko{\fe}$ be a family of holomorphic maps on $\skakko{|z| \ge R >0}$ such that
\begin{align*}
\fe(z) & \ee \taue z +1 +O(1/z) \\
\longrightarrow 
f_0(z) & \ee z +1 +O(1/z)
\end{align*}
uniformly as $\e \to 0$ where $\taue=1 + \e$. If $R \gg 0$, then for any $\e \in [0,1]$ there exists a holomorphic map $\ue:E_R \to \Cbar$ such that
$$
\ue(\fe(z)) \ee \taue \ue(z)+1
$$
and $\ue \to u_0$ uniformly on compact sets of $E_R$.
\end{thm}

Indeed, Ueda's original theorem in \cite{Ue} claims that a similar holds for any radial convergence $\taue \to 1$ outside the unit disk. In \cite{Ka5} an alternative proof is given and the error term $O(1/z)$ is refined to be $O(z^{-1/n})$ for any $n \ge 1$.

\begin{lem}\label{lem_u_0}
$u_0(z)=z(1+o(1))$ as $\rp z \to \infty$. 
\end{lem}
Indeed, it is well-known that if $f_0(z)=z+1+a_0/z+\cdots$ then the Fatou coordinate is of the form $u_0(z)=z-a_0 \log z+O(1)$. See \cite{Sh} for example.

\subsection{Small invariant paths joining perturbed periodic points}
For a degeneration pair $(f \to g)$ in Case (a) ($q=q'$), we may consider that the parabolic cycle $O_g$ is perturbed into the attracting cycle $O_f$ with the same period $l=l'$. (See \propref{prop_case_a_case_b}). In this case, the convergence $f^{lq} \to g^{lq}$ is viewed as 
$$
f^{lq}(z)=r^q z +z^q+O(z^{2q}) ~\longrightarrow~ g^{lq}(z)=z +z^q+O(z^{2q})
$$
with $r^q \nearrow 1$ through suitable local coordinates near $\beta_0 \in O_g$ as in Appendix A.2. 

By taking an additional linear coordinate change by $z \mapsto z/r$, we consider a family of holomorphic maps $\skakko{\fe}$ of the form
$$
\fe(z)=\lame z(1+z^q+O(z^{2q}))
$$
instead, where we set $r^q=\lame=1-\e \nearrow 1$. Then the local solution of $\fe(z)=z$ is $z=0$ or $z^q=\e+O(\e^2)$. The latter means $q$ symmetrically arrayed repelling fixed points are generated by the perturbation of a parabolic point with multiplicity $q+1$. Here we claim:

\begin{lem}\label{lem_small_inv_path}
For $\e \ll 1$, there exist $q$ $\fe$-invariant paths of diameter $O(\e^{1/q})$ joining the central attracting point $z=0$ and each of symmetrically arrayed repelling fixed points.
\end{lem}

\begin{pf}
First we show that $D:=\skakko{z: |z|^q \le \e/2}$ satisfies $\fe(D) \subset D^\cc$. By checking the real part of $\log \fe(z)$, we have 
$$
|\fe(z)| \ee \lame |z|(1+\rp z^q+O(z^{2q})).
$$
Since $\rp z^q \le \e/2$ on $D$, we have $|\fe(z)|=|z|(1-\e/2+O(\e^2))<|z|$.

Next we set
$$
E \dee \skakko{z: \frac{\e}{2} \le |z^q| \le 4\e 
~~\text{and}~~|\arg z^q| \le \frac{\pi}{3}}.
$$
Note that $E$ has $q$ connected components around the repelling fixed points. Now we claim that $E$ satisfies $\fe^{-1}(E) \subset E^\cc$. Since $\fe^{-1}$ is univalent near $0$, it is enough to show that $\fe^{-1}(\partial E) \subset E^\cc$. Set
\begin{align*}
e_1 &\dee \skakko{z: |z^q| = \frac{\e}{2} 
~~\text{and}~~|\arg z^q| \le \frac{\pi}{3}} \\
e_2 &\dee \skakko{z: |z^q| = 4\e 
~~\text{and}~~|\arg z^q| \le \frac{\pi}{3}} \\
e_3^\pm &\dee \skakko{z: \frac{\e}{2} \le |z^q| \le 4\e
~~\text{and}~~\arg z^q =\pm \frac{\pi}{3}}. 
\end{align*}
By checking $\log \fe^{-1}(z)$, we have 
\begin{align*}
|\fe^{-1}(z)| &\ee \lame^{-1} |z|(1- \lame^{-q}\rp z^q+O(z^{2q})) \\
\arg \fe^{-1}(z) &\ee \arg z -\lame^{-q}\ip z^q+O(z^{2q}).
\end{align*}
If $z \in e_1$, we have $\rp z^q \le \e/2$ thus $|\fe(z)|\ge |z|(1+\e/2+O(\e^2))>|z|$. If $z \in e_2$, we have $\rp z^q \ge 2\e$ thus $|\fe(z)|\le|z|(1-\e+O(\e^2))<|z|$. For $z \in e_3^\pm$, set $|z^q|=\rho$ with $\e/2 \le  \rho \le 4\e$. Then $
\arg \fe^{-1}(z) = \arg z \mp (\sqrt{3}/2)\rho(1 + O(\rho)).
$
Thus we have $\fe^{-1}(\partial E) \subset E^\cc$ in total.

Take any $q$ points $\skakko{z_1, \cdots, z_q}$ from each connected component of $e_1$. Let $\eta_j$ be the segment joining $z_j$ and $\fe(z_j)$. Then the path $\bigcup_{k \in \Z}\fe^{k}(\eta_j)$ has the desired property.
\QED
\end{pf}

\parag{Remark.}
In Case (b) ($q=1<q'$), the cycle $O_f'$ in Appendix A.2 is repelling. By taking $f^{-l'q'} \to g^{-l'q'}$ near $O_g$, we have a similar form of the convergence
$$
\fe(z)=\lame z(1+z^{q'}+O(z^{2q'}))
$$
to the case of $q=q'$, where $\lame=1-\e +O(\e^2) \in \Cstar$. This $\lame$ comes from the fact that the non-zero solutions of $\fe(z)=z$ has derivative $0<r<1$ (Since they are actually points in $O_f$ in a different coordinate.) One can easily check that the argument of \lemref{lem_small_inv_path} above also works for this $\fe$ and the statement is also true by replacing $q$ with $q'$.

\end{document}